\let\latexarabic\arabic
\let\latexdocument\document
\let\latexenddocument\enddocument

\RequirePackage[thmmarks]{ntheorem}
\makeatletter
\renewtheoremstyle{plain} 
  {\item[\hskip\labelsep \theorem@headerfont ##1\ \textup{##2}\theorem@separator]} 
  {\item[\hskip\labelsep \theorem@headerfont ##1\ \textup{##2}\ (##3)\theorem@separator]}
\makeatother

\documentclass[
  supplementary,
  lineno
]{biometrika}

\let\document\latexdocument
\let\enddocument\latexenddocument
\AtEndDocument{\printhistory}
\let\arabic\latexarabic
\def\rm{}
\makeatletter
\newtheoremstyle{MyNonumberplain}%
  {\item[\theorem@headerfont\hskip\labelsep ##1\theorem@separator]}%
  {\item[\theorem@headerfont\hskip\labelsep ##3\theorem@separator]}
\makeatother
\theoremstyle{MyNonumberplain}
\theorembodyfont{\upshape}
\newtheorem{Proof}{Proof}

\usepackage{amsfonts,amssymb}
\usepackage{bm}
\usepackage{mathtools}
\usepackage[center]{titlesec}
\usepackage{algorithm}
\usepackage{subcaption}
\usepackage[colorlinks=true,linkcolor=blue, citecolor=blue]{hyperref}

\usepackage{amsmath}

\usepackage{times}
\usepackage{natbib}

\numberwithin{equation}{section}

\newcommand{\rO}{\mathrm{O}}

\DeclareMathOperator*{\argmin}{argmin}

\theoremstyle{plain}
\newtheorem{thm}{Theorem}[section]
\newtheorem{cor}[thm]{Corollary}

\newtheorem{assu}[thm]{Assumption}

\newtheorem{defn}[thm]{Definition}
\newtheorem{exam}[thm]{Example}
\theoremstyle{remark}
\newtheorem{rem}[thm]{Remark}

\begin{document}
\nolinenumbers
\jname{}
\jyear{}
\jvol{}
\jnum{}


\markboth{X. Ding and Z. Zhou}{PACF for locally stationary time series}

\title{On the partial autocorrelation function for locally stationary time series: characterization, estimation and inference}

\author{Xiucai Ding}
\affil{Department of Statistics, University of California, Davis, Davis 95616 USA \email{xcading@ucdavis.edu}}

\author{\and Zhou Zhou}
\affil{Department of Statistical Sciences, University of Toronto, Toronto M5G 1X6, Canada \email{zhou.zhou@utoronto.ca}}

\maketitle

\begin{abstract}
For stationary time series, it is common to use the plots of partial autocorrelation function (PACF) or PACF-based tests to explore the temporal dependence structure of such processes. To our best knowledge, such analogs for non-stationary time series have not been fully established yet. In this paper, we fill this gap for locally stationary time series with short-range dependence. {First, we characterize the PACF locally in the time domain and show that the $j$th PACF, denoted as $\rho_{j}(t),$ decays with $j$ whose rate is adaptive to the temporal dependence of the time series $\{x_{i,n}\}$. Second, at time $i,$ we justify that the PACF $\rho_j(i/n)$ can be efficiently approximated by the best linear prediction coefficients via the Yule-Walker's equations. This allows us to study the PACF via ordinary least squares (OLS) locally. Third, we show that the PACF is smooth in time for locally stationary time series. We use the sieve method with OLS to estimate $\rho_j(\cdot)$} and construct some statistics to test the PACFs and infer the structures of the time series. These tests generalize and modify those used for stationary time series in \cite{BD}. Finally, a multiplier bootstrap algorithm is proposed for practical implementation and an $\mathtt R$ package $\mathtt {Sie2nts}$ is provided to implement our algorithm. Numerical simulations and real data analysis also confirm usefulness of our results. 
\end{abstract}


\vspace*{-1.5cm}
\section{Introduction}
The partial autocorrelation function (PACF) is one of the most popular and powerful tools for stationary time series modelling and analysis \cite{BD}. However, in the era of big data, as increasingly
longer time series are being collected, it has become more
appropriate to model many of those series as locally stationary processes whose data generating mechanisms evolve smoothly over time. In this setting, the effectiveness of the classical PACF deteriorates and it is of urgent demand to establish the theories of PACF for locally stationary time series. 

Even though there exists a rich body of literature on locally stationary time series analysis, see \cite{dahlhaus2012locally,dahlhaus2019towards} for a review, much less has been studied related to the PACF.  In \cite{DEGERINE200346}, the authors generalized the characterization of PACFs via some useful decomposition as introduced in \cite{ramsey1974characterization} to general non-stationary processes. They also briefly discussed how to estimate the PACFs based on a generalized Levinson-Durbin algorithm when the autocovariance function is given. More recently, in \cite{10.1214/20-EJS1748}, by generalizing the partial autocorrelations of stationary processes to locally stationary time series from a wavelet spectrum perspective, the authors provided two new estimators for the local PACFs. The consistency of the wavelet-based estimator and the asymptotic distribution of the windowed estimator under Gaussian assumption have also been studied.  However, the decay speed of local PACFs as a function of the lags has not been established and a direct time-domain characterization of the PACFs of locally stationary processes has not been fully investigated. Moreover, the inference for PACFs of locally stationary time series, for example significance tests and PACF-based Portmanteau tests,  are still missing in the literature.

Motivated by the above challenges, in this paper, we aim to systematically study the theories of PACFs for locally stationary time series. For characterization, in contrast to \cite{DEGERINE200346,10.1214/20-EJS1748}, we define PACFs  for general locally stationary time series {in the time domain using stationary approximations at each time point (c.f. Definition \ref{defn_locallystationary}). There are several advantages in using this characterization. First, since the time series is approximately stationary locally, the lower-order PACFs can be well approximated by the best short-term linear prediction coefficients via the Yule-Walker equations; see (\ref{eq_functionfinal}) and (\ref{eq_smoothapproximation}). This connection not only allows us to study the PACF via ordinary least squares (OLS) estimation but also enables us to establish the decay properties of the PACFs which are adaptive to the temporal dependence decay of the time series; see (\ref{eq_decayrate}).}
Second, the smoothness of the locally stationary covariance structure can be easily translated  to that of the PACFs; see (\ref{eq_smoothapproximation}). Therefore, it suffices to estimate some smooth functions $\rho_j(\cdot)$ at different lags $j.$ More specifically, {together with the OLS,} the smooth PACFs can be estimated adaptively using the nonparametric method of sieves via flexible choices of basis functions such as the wavelets and the orthogonal polynomials \cite{chen2007large}; see Section \ref{sec_estimation}. Theoretically, {under mild assumptions, the estimators are  consistent uniformly in the time domain (c.f. Theorem \ref{thm_consistency}).} Third, based on the time-domain characterization and the OLS form of the sieve estimators, one can further conduct various tests on the PACFs.
%
%
For example, one can perform a  white noise Portmanteau test or significance tests on some PACFs (e.g., checking the order for an AR process) uniformly over time and lags. Both tests have not been fully studied yet under the locally stationary time series framework.  We establish the asymptotic normality and conduct a power analysis for the tests (c.f. Theorems \ref{lem_controltbtrho0} and \ref{lem_controltbtrho}).  We also propose a multiplier bootstrap procedure for practical implementation (c.f. Algorithm \ref{alg:boostrapping}). Numerical simulations and real data analysis are provided to illustrate the usefulness of our results and an $\mathtt{R}$ package $\mathtt{Sie2nts}$ is provided for users. Since our method covers stationary time series as a special case, we promote the use of the $\mathtt{Sie2nts}$ package instead of the default $\mathtt{pacf}$ function in  $\mathtt{R}$ which only handles stationary time series.

The paper is organized as follows. In Section \ref{sec_characterization}, we provide the characterization of the PACFs for locally stationary time series and study their asymptotic properties. In Section \ref{sec_estimationinference}, we introduce our estimator for PACFs based on the nonparametric sieve method and inference procedures based on multiplier bootstrap. In Section \ref{sec_theorecticalanalysis}, we provide theoretical analysis for our estimation and inference procedures. Numerical simulations and real data analysis are offered in Section \ref{sec_simulation}.  More details are provided in our online supplement \cite{supplement}. Especially, technical proofs are deferred to Section \ref{sec_techinicalproof}, further discussions and remarks are provided in Section \ref{sec_additionalremark}, some tuning parameter selection algorithm is listed in Section \ref{sec:choiceparameter} and additional simulation results are offered in Section \ref{sec_additionalsimulationresults}.  


\vspace{5pt}

\noindent{\bf Conventions.} For a random variable $x \in \mathbb{R}$ and some constant $q \geq 1,$ we denote by $\|x \|_q:=(\mathbb{E}|x|^q)^{1/q}$ its $L^q$ norm. We simply write $\| x\| \equiv \|x \|_2$ when $q=2.$ For any deterministic vector $\bm{x} =(x_1, x_2, \cdots, x_p)^* \in \mathbb{R}^p,$ we use $|x|=\sqrt{\sum_{i=1}^p x_i^2}$ for its $\ell_2$ (or Euclidean) norm. For any matrix $A,$ we use $\| A \|$ to stand for its operator norm. For two sequences of deterministic positive values $\{a_n\}$ and $\{b_n\}$, we write $a_n=\mathrm{O}(b_n)$ if $a_n \leq C b_n$ for some positive constant $C>0$. Moreover, we write $a_n=\mathrm{o}(b_n)$ if $a_n \leq c_n b_n$ for some positive sequence $c_n \downarrow 0.$ For a sequence of random variables $\{x_n\}$ and positive real values $\{a_n\},$ we use the notation $x_n=\rO_{\ell^q}(a_n)$ to state that  $x_n/a_n$ is bounded in $L^q$ norm; that is $\| x_n/a_n \|_q \leq C$ for some finite constant $C>0.$ If $q=1,$ we simply write $x_n=\mathrm{O}_{\mathbb{P}}(a_n).$ We use $\mathbf{I}_c$ for a $c \times c$ identity matrix. We use $C^d([0,1])$ for the function space on $[0,1]$ of continuous functions that have continuous first $d$ derivatives.

\vspace*{-0.5cm}
\section{Characterization of PACF for locally stationary time series}\label{sec_characterization}

In this section, we provide the characterization of PACF of locally stationary time series and study its properties. Suppose that we observe $\{x_{i,n}\}_{i=1}^n$. For simplicity and without loss of generality, we assume the time series is centered (i.e., mean is zero). Till the end of the paper, for notional simplicity,  we always write $x_i \equiv x_{i,n}.$ From time to time, we will emphasize the dependence on $n$ for various quantities.

In this paper, we focus on a general class of locally stationary time series following definition introduced in \cite{DZar}. It covers many commonly used locally stationary time series models in the literature. See for instance \cite{DZ1,dahlhaus2019towards,dahlhaus2012locally,RSP,KPF,DPV,dette2020prediction,MR3097614,WZ1}. We refer the readers to Example \ref{exam_timeseries} in our supplement for more details. 
\vspace*{4pt}

\begin{defn}[Locally stationary time series and its PACF]\label{defn_locallystationary} A non-stationary time series  $\{x_{i}\}$ is a locally stationary time series (in covariance) if there exists a function $\gamma(t, k): [0,1] \times \mathbb{N} \rightarrow \mathbb{R}$ such that \vspace*{-0.1cm}
\begin{equation}\label{eq_covdefn}
\operatorname{Cov}(x_{i},x_{j})=\gamma(t_i, |i-j|)+\rO \left( \frac{|i-j|+1}{n} \right), \ t_i=\frac{i}{n}.
\end{equation} 
Moreover, we assume that $\gamma$ is Lipschitz continuous in $t$ and for any fixed $t \in[0,1],$ $\gamma(t,\cdot)$ is the autocovariance function (ACF) of some stationary process {whose $j$th order PACF is denoted as $\rho_j(t)$. We shall also call $\rho_j(t)$ the $j$th order PACF at rescaled time $t$ of $\{x_i\}.$}
\end{defn}

Observe that Definition  \ref{defn_locallystationary} essentially means that the covariance structure of $\{x_{i}\}$ can be well approximated locally by that of a stationary process. Consequently, at each rescaled time $t$, the PACF of $\{x_{i}\}$ is defined by that of the approximating stationary process.  Now we study the PACF introduced above.  Since only one realization of the series $\{x_{i}\}$ is available, we will need to assume short-range dependence that for some $\tau>1,$
\begin{equation}\label{eq_polynomialdecay}
\max_{k,n}\left| \operatorname{Cov}(x_{k,n}, x_{k+r,n}) \right| =\mathrm{O} \left( |r|^{-\tau} \right) \ \text{and} \ \sup_t|\gamma(t, |r|)|=\mathrm{O}(|r|^{-\tau}),
\end{equation}
and add additional regularity conditions. These will be summarized in Assumption \ref{assum_mainassumption} after some necessary notations are introduced. For the autocovariance function $\gamma(t, \cdot)$ in (\ref{eq_covdefn}), given a lag $j,$ we define a vector of functions $\bm{\phi}_j(t)=(\bm{\phi}_{j,1}(t), \cdots, \bm{\phi}_{j,j}(t))^*: [0,1] \rightarrow \mathbb{R}^j$ via the local Yule-Walker equation 
\begin{equation}\label{eq_defnphjt}
\bm{\phi}_j(t)=\Gamma_j(t)^{-1} \bm{\nu}_j(t),
\end{equation}
where $\Gamma_j(t) \in \mathbb{R}^{j \times j}$ is a symmetric matrix whose $(k,l)$ entry is defined as $\Gamma_j^{(k,l)}(t):=\gamma(t, |k-l|),$ and $\bm{\nu}_j(t) \in \mathbb{R}^j$ is vector whose $k$th entry is defined as $\bm{\nu}_j^{(k)}(t):=\gamma(t, k).$ {Consequently, we can write \cite{BD} }
\begin{equation}\label{eq_functionfinal}
\rho_j(t)=\bm{\phi}_{j,j}(t). 
\end{equation}

{Note that for stationary time series the PACF is closely related to the best linear forecast coefficients of the process (\cite{BD}).  The representation (\ref{eq_functionfinal}) together with  the smoothness of the time series covariance structure imply that we can study the PACF of locally stationary time series via local best linear forecasts and therefore a simple local regression analysis.} Specifically, for all $1 \leq j \leq i-1,$ denote the $j$th order best linear forecast of  $x_{i}$ as 
$\widehat{x}_{i,j}=\sum_{k=1}^j \phi_{ik,j} x_{i-k},$
where $\phi_{ik,j} \equiv \phi_{ik,j,n}, \ 1 \leq k \leq j,$ are the best linear forecast coefficients. Define the residual as $\epsilon_{i,j} \equiv \epsilon_{i,j,n}:=x_i-\widehat{x}_{i,j}.$ 
We now write \vspace*{-0.1cm}
\begin{equation}\label{eq_jorderlocal}
x_i=\sum_{k=1}^j \phi_{ik,j} x_{i-k}+\epsilon_{i,j}.   
\end{equation}
Using (\ref{eq_jorderlocal}), we introduce the following notations \vspace*{-0.1cm}
\begin{equation}\label{eq_rijdefinition}
\rho_{i,j}=\phi_{ij,j}, \ 1 \leq j \leq i-1. 
\end{equation}
{The following theorem studies the uniform decay property of the PACF and builds the connection between the PACF defined in Definition \ref{defn_locallystationary} and $\rho_{i,j}$ in (\ref{eq_rijdefinition}) for the locally stationary time series.}

\vspace*{4pt}

\begin{thm}\label{lem_pacf}
For the locally stationary time series $\{x_i\}$ satisfying Definition \ref{defn_locallystationary}, suppose Assumption \ref{assum_mainassumption} holds. Then we have the followings holds. 
\begin{enumerate}
\item For $\tau$ in (\ref{eq_polynomialdecay}), we have that  for $\rho_j(t)$ in Definition \ref{defn_locallystationary} and $\rho_{i,j}$ in (\ref{eq_rijdefinition}) \vspace*{-0.1cm}
\begin{equation}\label{eq_decayrate}
\sup_{i>j} |\rho_{i,j}|=\rO\left( \left(\frac{\log j+1}{j} \right)^{\tau-1} \right) \ \text{and} \ \sup_t |\rho_j(t)|=\rO\left( \left(\frac{\log j+1}{j} \right)^{\tau-1} \right).
\end{equation}
\item $\rho_j(t) \in C^d([0,1]),$ {for some integer $d>0$ defined in (\ref{eq_definitiond})} and \vspace*{-0.1cm}
\begin{equation}\label{eq_smoothapproximation}
\sup_{i>j} \left| \rho_j \left(\frac{i}{n} \right)-\rho_{i,j} \right|=\rO \left( \frac{j^2}{n} \right).
\end{equation}
\end{enumerate} 
%
%
\end{thm}

\begin{rem}\label{rem_key}
Several remarks on Theorem \ref{lem_pacf} are in order. First, (\ref{eq_decayrate})  implies that, uniformly over time, the PACF decays polynomially fast to 0 as a function of the lag with the speed adaptive to that of the autocovariance.  If $\tau$ is sufficiently large, this implies that practically we only need to consider the first few lags of the PACF. In practice, people are usually concerned with the cutoff that {$ \sup_{t} |\rho_{j}(t)| \gg n^{-1/2}.$} In this sense, we conclude from (\ref{eq_decayrate}) that we only need to focus on the lags for $j=\rO \left( n^{\frac{1}{2(\tau-1)}} \right).$ Second, (\ref{eq_smoothapproximation}) demonstrates that under Definition \ref{defn_locallystationary} and Assumption \ref{assum_mainassumption}, the PACFs can be well approximated by $\rho_{i,j}$ if $j$ is not too large. Observe that $\rho_{i,j}$ is defined by local best linear forecasts and therefore is closely related to OLS. Therefore, (\ref{eq_smoothapproximation}) lays the theoretical foundation for the estimation and inference procedures in Section \ref{sec_estimationinference}. 
%
%
\end{rem}

\vspace*{-0.5cm}
\section{Estimation and inference procedures}\label{sec_estimationinference}
In this section, we provide the procedures of the estimation and inference of the PACFs of locally stationary time series. The theoretical justifications will be provided in Section \ref{sec_theorecticalanalysis}. As mentioned in Remark \ref{rem_key}, in what follows, we shall only consider lag $j \leq j_*$ with
\begin{equation}\label{j_assumption}
j_* \asymp  n^{\frac{1}{2(\tau-1)}} . 
\end{equation}

\subsection{Sieve nonparametric estimation}\label{sec_estimation}
In this section, we estimate the PACFs $\rho_j(\cdot).$ As proved in Theorem \ref{lem_pacf}, since $\rho_j(t) \in C^d([0,1]),$ it is natural for us to approximate it via a finite diverging term basis expansion using the method of sieves as in \cite{chen2007large, DZ1, DZar}. Recall $\rho_j(t)$ is defined via $\bm{\phi}_j(t)$ in (\ref{eq_functionfinal}). Now we will work with $\bm{\phi}_j(t).$ 

According to \cite[Section 2.3]{chen2007large}, we have that for some pre-chosen orthonormal basis functions on $[0,1]$, denoted as $\{\alpha_k(t)\}$ \vspace*{-0.1cm}
\begin{equation}\label{eq_rhocxh}
\bm{\phi}_{j,l}(t):=\sum_{k=1}^c a_{jk,l} \alpha_k(t)+\rO(c^{-d}), \ 1 \leq l\leq j,
\end{equation}
where $c$ is the number of basis functions. Here $\{a_{jk,l}\}$ are the coefficients to be estimated. In fact, using \cite[Theorem 2.11]{DZar} and a discussion similar to (3.8) therein, we have by \eqref{eq_jorderlocal} that \vspace*{-0.1cm} 
\begin{equation}\label{eq_representationxi}
x_i=\sum_{l=1}^j \sum_{k=1}^c a_{jk,l} z_{kl}+\epsilon_{i,j}+\rO_{\ell^2}\left(j^{2}/n+jc^{-d}\right), \ i>j,
\end{equation}
where $\epsilon_{i,j}$ is defined in (\ref{eq_jorderlocal}) and $z_{kl}(i/n)=\alpha_k(i/n)x_{i-l}$. Observe that  a key component in \eqref{eq_representationxi} is the approximation of $\rho_j(\cdot)$ by local best linear forecast coefficients as in \eqref{eq_smoothapproximation}.

Using (\ref{eq_representationxi}), we can estimate all the coefficients $a_{jk,l}$'s using one OLS regression. In particular, we stack $a_{jk,l},  1 \leq l \leq j, 1 \leq k \leq c$ as a vector $\bm{\beta}_j \in \mathbb{R}^{jc},$ then the OLS estimator for $\bm{\beta}_j$ can be written as $\widehat{\bm{\beta}}_j:=(Y_j^* Y_j)^{-1} Y_j^* \bm{x}_j,$ where $Y_j^* \in \mathbb{R}^{jc \times (n-j)}$ is the design matrix of (\ref{eq_representationxi}) and $\bm{x}_j=(x_{j+1}, \cdots, x_n)^* \in \mathbb{R}^{n-j}.$ 

 After estimating the $a_{jk,l}$'s, $\rho_j(t) \equiv \bm{\phi}_{j,j}(t)$ is estimated using (\ref{eq_rhocxh}) as
 \begin{equation}\label{eq_pacfestimate}
 \widehat{\rho}_j(t)=\widehat{\bm{\beta}}_j^* \mathbb{B}_{j,j}(t),
 \end{equation}
 where $\mathbb{B}_{j,l}(t) \in \mathbb{R}^{jc}$ has $j$ blocks and the $l$th block is $\mathbf{B}(t):=(\alpha_1(t), \cdots, \alpha_c(t))^* \in \mathbb{R}^c, 1 \leq l \leq j$ and zeros otherwise.
 
\subsection{Multiplier bootstrap based inference} 
 
In this subsection we propose a multiplier bootstrap procedure to infer the PACFs. Statistical inference of the PACFs plays an important role in stationary time series analysis. For example, it can be used to determine the order of an AR process and check whether the time series (or residuals after an ARIMA model fitting) is white noise. We refer the readers to Chapter 3 of \cite{shumway2017time} for more details. However, the analogs for locally stationary time series are largely missing. We aim to fill the gap in this section.

Based on our estimators in (\ref{eq_pacfestimate}), we can conduct various important tests on $\rho_j(t)$ in (\ref{eq_functionfinal}). For example, we can test whether the PACFs are identical to functions of interest that $\rho_j(t) = f_j(t)$ for some given functions $f_j(t).$  For another instance, we can check whether a group of  the PACFs are time-invariant that $\rho_j(t) \equiv \rho_j, \ h_1 \leq j \leq h_2$ for some integers $h_1$ and $h_2,$ where $ \rho_j(t) \equiv \rho_j$ means $ \rho_j(t)=\rho_j$ for all $t$. Note that if we set $\rho_j=0,$ it reduces to the testing the significance of the PACFs. Especially, when $h_1=1$ and $h_2=\infty,$ it is equivalent to test whether the time series is white noise. While we are able to conduct several different tests on the PACFs, in this paper, motivated by their applications in model selection and goodness of fit,  for conciseness, we will focus on two important such tests.

First, we are interested in testing 
\begin{equation}\label{eq_firstnull}
\mathbf{H}_{01}: \ \rho_j(t) \equiv 0 \ \  \text{vs} \ \  \mathbf{H}_{a1}: \ \rho_j(t) \not \equiv 0, 
\ \text{for some} \ j \geq 1.
\end{equation}
 The hypothesis (\ref{eq_firstnull}) tests the significance of a single PACF. Similar to the stationary setting in \cite{BD}, it can be used to select the order of a locally stationary AR process; see Remark \ref{rem_app_determinep} of our supplement for more discussions.

Second, we are also interested in testing the significance for all the lags that 
 \begin{equation}\label{eq_testnotwo}
 \mathbf{H}_{02}: \rho_{j}(t)  \equiv 0 \ \text{for all} \ j \geq 1 \ \ \text{vs} \ \ \mathbf{H}_{a2}: \rho_k(t) \not \equiv 0 \ \text{for some} \ k \geq 1. 
 \end{equation}
The hypothesis (\ref{eq_testnotwo}) tests white noise (or lack
of serial correlation) of the underlying locally time series. For stationary white noise, the well-known Box–Pierce (BP) test statistic \cite{box1970distribution} with fixed lag truncation number is probably the most commonly used statistic. Later on, such a test was extended to locally stationary white noise in \cite{goerg2012testing}. We emphasize that the Portmanteau-type BP tests in \cite{box1970distribution,goerg2012testing} used the autocorrelation functions (ACFs) instead of the PACFs. However, the estimation of ACFs for locally stationary time series involves the estimation of the time-varying marginal variances which requires the choice of additional tuning parameters and may lead to deteriorated estimation accuracy in finite samples. {Inspired by the above challenges and the discussions of Section 9.4 of \cite{BD}, we will propose a PACF-based Portmanteau test.}

We mention again that as discussed in Remark \ref{rem_key}, when $\tau$ is large, we have that $\sup_t\rho_j(t)=\mathrm{o}(n^{-1/2})$ for $j \gg j^*$ defined in (\ref{j_assumption}). Therefore, from an inferential viewpoint, $\rho_j(\cdot)$ for $j \gg j^*$ can be effectively treated as zero. Consequently, we only need to consider the setting $j \leq j^*$ for (\ref{eq_firstnull}) and (\ref{eq_testnotwo}) once $j^*$ is identified. 
\vspace*{-0.3cm}  
\subsubsection{Test statistics}
In this section,  we propose the test statistics for (\ref{eq_firstnull}) and (\ref{eq_testnotwo}). First,  when the null hypothesis in  (\ref{eq_firstnull}) holds, the following statistic $T_1 \equiv T_{1}(j)$ should be small \vspace*{-0.1cm}
\begin{equation}\label{eq_teststatisticone}
T_1 \equiv T_{1}(j):=\int_0^1 \widehat{\rho}_j^2(t) \mathrm{d} t. 
\end{equation} 
Therefore, it is natural to use (\ref{eq_teststatisticone})  to test (\ref{eq_firstnull}).

 Second, to test (\ref{eq_testnotwo}), motivated by the BP test \cite{box1970distribution}, we may want to directly use the following statistic \vspace*{-0.1cm}
\begin{equation}\label{eq_Tbp}
T_{\text{BP}} \equiv T_{\text{BP}}(\mathsf{h}):=\sum_{k=1}^{\mathsf{h}} \int_0^1  \widehat{\rho}_k(t)^2 \mathrm{d} t, \ \text{for some large} \ \mathsf{h} \geq j^*,
\end{equation}
where we recall $j^*$ in (\ref{j_assumption}). Even though it is natural to use $T_{\text{BP}}$, as described in Section \ref{sec_estimation}, in order to obtain its value, we need to do $\mathsf{h}$ high dimensional OLS regressions which can be computationally expensively, especially when $\mathsf{h}$ is large (or equivalently, $\tau$ is small). To address this issue, we consider an $\mathsf{h}$ order best linear prediction  as in (\ref{eq_jorderlocal}). That is, $x_i=\sum_{k=1}^{\mathsf{h}} \phi_{ik,\mathsf{h}} x_{i-k}+\epsilon_{i,\mathsf{h}}. $ According to Theorem 2.11 of \cite{DZar}, by setting $j=\mathsf{h}$ in (\ref{eq_defnphjt}), we find that \vspace*{-0.1cm}
\begin{equation}\label{eq_smoothapproximationar}
x_i= \sum_{k=1}^{\mathsf{h}} \phi_{\mathsf{h},k}(i/n) x_{i-k}+\epsilon_{i, \mathsf{h}}+\rO_{\ell^2}\left( \mathsf{h}^{2}/n\right). 
\end{equation}
The smooth coefficients $\{\phi_{\mathsf{h},k}(\cdot)\}$ can be estimated via the sieve method using only one high-dimensional OLS as in Section \ref{sec_estimation} whose estimators are denoted as $\{\widehat{\phi}_{\mathsf{h},k}(\cdot)\}.$

Now we define another statistic \vspace*{-0.1cm}
\begin{equation}\label{eq_defnt2}
T_2 \equiv T_2(\mathsf{h}):=\sum_{k=1}^{\mathsf{h}} \int_0^1 \widehat{\phi}_{\mathsf{h},k}^2(t) \mathrm{d} t. 
\end{equation} 
As will be seen later in Theorem \ref{lem_controltbtrho} below, when (\ref{eq_testnotwo}) holds, under some mild conditions on $\mathsf{h}$ {(c.f. (\ref{eq_anotherbound}))}, $T_2$ will be close to $T_{\text{BP}}$ while the calculation of $T_2$ only needs one OLS regression. Therefore, we will use (\ref{eq_defnt2}) to test (\ref{eq_testnotwo}).  
\vspace*{-0.3cm}
\subsubsection{Practical implementation}
We point out that it is still difficult to directly use $T_{1}$ and $T_2$ since the variances in their limiting Gaussian distributions are usually hard to estimate and hence  plug-in estimators are unavailable; see Theorems \ref{lem_controltbtrho0} and \ref{lem_controltbtrho} for more details. To address this issue, we utilize the multiplier bootstrap procedure as in \cite{DZar,ZZ1}. We first explain how to construct the bootstrapped statistics. Using the sieve estimates as in Section \ref{sec_estimation}, denote the residual \vspace*{-0.1cm}
\begin{equation}\label{eq_residualestimate}
\widehat{\epsilon}_{i,\ell}:=x_i-\sum_{k=1}^{\ell} \widehat{\phi}_k(i/n) x_{i-k}, 
\end{equation}   
where $\ell=j$ for $T_1$ and $\ell=\mathsf{h}$ for $T_2.$ 

Let  $\widehat{\bm{w}}_i=\bm{x}_{\ell,i} \widehat{\epsilon}_{i,\ell},$ where $\bm{x}_{\ell,i}=(x_{i-1}, \cdots, x_{i-\ell})^*.$ Given a block size $m,$ we denote $\widehat{\Phi} \equiv \widehat{\Phi}(\ell, m) \in \mathbb{R}^{\ell c}$ as \vspace*{-0.1cm} 
\begin{equation}\label{eq_sampledstatistic}
\widehat{\Phi}:=\frac{1}{\sqrt{n-m-\ell+1} \sqrt{m}} \sum_{i=\ell+1}^{n-m} \left[ \left( \sum_{j=i}^{i+m} \widehat{\bm{w}}_{\ell,j}\right) \otimes \left( \mathbf{B}(i/n) \right) \right] R_i,
\end{equation}
{where $\otimes$ is the Kronecker product and} $R_i, \ell+1 \leq i \leq n-m$ are i.i.d. standard Gaussian random variables which are independent of the observed time series. Recall the discussions around (\ref{eq_pacfestimate}) and $Y_{\ell}^*$ is the design matrix. Denote $\widehat{\Sigma}_1:=\frac{1}{n} Y_{j}^* Y_{j}$ and $\widehat{\Sigma}_2:=\frac{1}{n} Y_{\mathsf{h}}^* Y_{\mathsf{h}}.$ Moreover, let  $\mathbf{M} \in \mathbb{R}^{\ell c \times \ell c}$ be  a diagonal block matrix whose only non-zero part is the identity matrix lies in the last diagonal block. Inspired by Remark \ref{rem_keykeykey} below, we use the following statistics $\widehat{\mathcal{T}}_k, k=1,2,$ to mimic the distribution of $n T_k, k=1,2.$ Let $\widehat{\Phi}_1$ be constructed as in (\ref{eq_sampledstatistic}) using $\ell=j$ and $\widehat{\Phi}_2$ be constructed  using $\ell=\mathsf{h}.$ Then we denote
\begin{equation}\label{eq_finalfinalsamplesamplestatistic}
\widehat{\mathcal{T}}_1=\widehat{\Phi}_1^* \widehat{\Sigma}^{-1} \mathbf{M} \widehat{\Sigma}_1^{-1} \widehat{\Phi}_1,  \ \widehat{\mathcal{T}}_2=\widehat{\Phi}_2^* \widehat{\Sigma}_2^{-2}  \widehat{\Phi}_2.  
\end{equation}
{We point out that for the implementation, one needs to select some large value of $\mathsf{h}$ to construct $T_2.$ We discuss how to choose this parameter in Section \ref{sec:choiceparameter} of our supplement.}

Finally, based on the above result, we propose the following Algorithm \ref{alg:boostrapping} for the practical implementation. Note that in order to implement Algorithm \ref{alg:boostrapping}, two tuning parameters, the number of basis functions $c$ and the block size $m$, have to been chosen properly. {In our $\mathtt{R}$ package $\mathtt{Sie2nts}$, these parameters can be chosen automatically using the function $\mathtt{auto.pacf.test}$ according to the methods provided in Section \ref{sec:choiceparameter} of our supplement \cite{supplement}.}

\begin{algorithm}[!ht]
\caption{\bf Multiplier  Bootstrap}
\label{alg:boostrapping}

\normalsize
\begin{flushleft}
\noindent{\bf Inputs:} The lag $j$ for (\ref{eq_teststatisticone}) or $\mathsf{h}$ for (\ref{eq_defnt2}), type I error rate $\alpha,$ tuning parameters $c$ and $m$ chosen by the data-driven procedure demonstrated in Section \ref{sec:choiceparameter} of our supplement, time series $\{x_i\},$ and sieve basis functions.

\noindent{\bf Step one:} Compute $\widehat{\Sigma}_1^{-1}$ using $n(Y_{j}^*Y_{j})^{-1}$ for (\ref{eq_firstnull}) or $\widehat{\Sigma}_2^{-1}$ using $n(Y_{\mathsf{h}}^*Y_{\mathsf{h}})^{-1}$ for (\ref{eq_firstnull}),  and the residuals $\{\widehat{\epsilon}_{i, \ell}\}_{i=\ell+1}^n$ according to (\ref{eq_residualestimate}) with $\ell=j$ for (\ref{eq_teststatisticone}) or $\ell=\mathsf{h}$ for (\ref{eq_defnt2}). 

\noindent{\bf Step two:}  Generate $B$ (say 1,000) i.i.d. copies of $\{\widehat{\Phi}^{(s)}\}_{s=1}^B$ according to (\ref{eq_sampledstatistic}). Compute $\widehat{\mathcal{T}}_k^s$, $k=1,2$,\ $s=1,2,\cdots, B,$  correspondingly as in (\ref{eq_finalfinalsamplesamplestatistic}). 

\noindent{\bf Step three:} Let $\widehat{\mathcal{T}}^{(1)}_k \leq \widehat{\mathcal{T}}^{(2)}_k \leq \cdots \leq \widehat{\mathcal{T}}^{(B)}_k$ be the order statistics of $\widehat{\mathcal{T}}^{s}_k, s=1,2,\cdots, B.$ Reject $\mathbf{H}_{01}$ in (\ref{eq_firstnull}) at the level $\alpha$ if $nT_1>\widehat{\mathcal{T}}_1^{(\lfloor B(1-\alpha)\rfloor)},$ where $\lfloor x \rfloor$ denotes the largest integer smaller or equal to $x.$  Reject $\mathbf{H}_{02}$ in (\ref{eq_testnotwo}) at the level $\alpha$ if $nT_2>\widehat{\mathcal{T}}_2^{(\lfloor B(1-\alpha)\rfloor)}.$ Let $B_k^*=\max\{r: \widehat{\mathcal{T}}^{r}_k \leq nT_k\}, k=1,2.$

\noindent{\bf Output:} $p$-value of the tests (\ref{eq_firstnull}) and (\ref{eq_testnotwo}) can be computed, respectively as $1-\frac{B^*_k}{B}, k=1,2.$
\end{flushleft}
\end{algorithm}

\vspace*{-0.5cm}
\section{Theoretical analysis}\label{sec_theorecticalanalysis}

In this section, we provide some theoretical analysis on our estimation and inference procedures. Till the end of the paper, for notational simplicity, we assume that the locally stationary time series admits the general physical representation equipped with the physical dependence measures (see (\ref{eq_physcialrepresentation}) and (\ref{eq_dependencemeasure}) of our supplement). In addition, we need the following assumptions. 
\vspace*{4pt}


\begin{assu}\label{assum_mainassumption}
Throughout the paper, we suppose the followings holds:
\begin{enumerate}
\item[(1).] For all sufficiently large $n \in \mathbb{N}, $ we assume that there exists a universal constant $\kappa>0$ that 
\begin{equation}\label{eq_defnkappa}
\lambda_n(\operatorname{Cov}(x_{1}, \cdots, x_{n})) \geq \kappa,
\end{equation} 
where $\lambda_n(\cdot)$ is the smallest eigenvalue of the given matrix and $\operatorname{Cov}(\cdot)$ is the covariance matrix of the given vector. 
\item[(2).] For all $n \in \mathbb{N},$ $1 \leq k \leq n$ and $-k+1 \leq r \leq n-k,$ we assume that there exists some constant $\tau>1$ such that (\ref{eq_polynomialdecay}) holds. In addition, we assume that $\sup_{i,n} \mathbb{E} |x_{i}|<\infty.$
\item[(3).] For some given integer $d>0$, we assume that for any $j \geq 0$
\begin{equation}\label{eq_definitiond}
\gamma(t,j) \in C^d([0,1]). 
\end{equation} 
\end{enumerate}
\end{assu}

The conditions in the above assumption are mild and can be satisfied by many commonly used locally stationary time series. 
Due to space constraint, we leave some discussions to Remark \ref{rem_assumptionremark} of supplement \cite{supplement}.
\subsection{Uniform consistency}

In what follows, we establish the consistency for our estimators. Recall $\mathbf{B}(t)$ below (\ref{eq_pacfestimate}). Denote \vspace*{-0.5cm}
 \begin{equation}\label{eq_parameteroneone}
 \xi_c=\sup_t \sup_{1 \leq i \leq c} |\alpha_i(t)|, \ \zeta_c=\sup_t |\mathbf{B}(t)|. 
\end{equation}  
The following mild assumption will be needed to ensure a consistent estimation, which has been used frequently in the literature, see \cite{DZ1, DZar,ding2021simultaneous,MR3097614}. Recall $\gamma(\cdot, \cdot)$ in (\ref{eq_covdefn}).
For all $j=1,2,\cdots, j_*,$ denote $\Sigma^{j}(t) \in \mathbb{R}^{j \times j}$ whose $(k,l)$th entry is $\Sigma_{kl}^{(j)}(t)=\gamma(t,|k-l|).$
\vspace*{4pt}
 
\begin{assu}\label{assu_integratedeigenvalueassumption}
For $j=1,2,3,\cdots, j_*,$ denote the long-run integrated covariance matrix as \vspace*{-0.1cm}
\begin{equation}\label{eq_longruncovariancematrixdefinition}
\Sigma^{(j)}=\int_0^1 \Sigma^{(j)}(t) \otimes \left( \mathbf{B}(t) \mathbf{B}(t)^* \right) \mathrm{d} t,
\end{equation}
where {we recall that} $\otimes$ is the Kronecker product.  We assume that the eigenvalues of $\Sigma^{(j)}$ are bounded above and also away from zero by some universal constants. 
\end{assu}
Then we proceed to state the main results of this section. Recall (\ref{eq_parameteroneone}). Denote \vspace*{-0.1cm}
\begin{equation}\label{eq_psiparameter}
\Psi(j,c):=j \xi_c \zeta_c \sqrt{\frac{c}{n}}  \left(1+ \frac{j^{2}}{n}+jc^{-d} \right)+c^{-d}.
\end{equation} 

\begin{thm} \label{thm_consistency}
Suppose Assumptions \ref{assum_mainassumption} and \ref{assu_integratedeigenvalueassumption}  hold true. Moreover, for $j \leq j^*$ satisfying \vspace*{-0.1cm} 
\begin{equation}\label{eq_rateassumption}
jc \left( \frac{\xi_c^2}{\sqrt{n}}+\frac{\xi_c^2 n^{\frac{2}{\tau+1}}}{n} \right)=\mathrm{o}(1),
\end{equation}
we have that for our estimator (\ref{eq_pacfestimate}) \vspace*{-0.1cm}
\begin{equation}\label{eq_consistentrate}
\sup_{t} \left|\widehat{\rho}_j(t)- \rho_j(t) \right|=\rO_{\mathbb{P}}\left( \Psi(j,c) \right). 
\end{equation}
\end{thm}

\begin{rem} \label{rem_consist}
Theorem \ref{thm_consistency} implies that our proposed  estimator (\ref{eq_pacfestimate}) is uniformly consistent under mild conditions. First, the condition (\ref{eq_rateassumption}) ensures that $n^{-1} Y_j^* Y_j$ in the OLS estimator will convergence to $\Sigma^{(j)}$ which guarantees the regular behavior of $\widehat{\bm{\beta}}_j.$ In fact, (\ref{eq_rateassumption}) can be easily satisfied.  Note that $\xi_c$ and $\zeta_c$ can be calculated for the specific sieve basis functions and so does the convergence rate in (\ref{eq_consistentrate}). For example, when $\{\alpha_k(t)\}$ are Fourier basis functions or normalized orthogonal polynomials, $\xi_c=\mathrm{O}(1)$ and $\zeta_c=\mathrm{O}(\sqrt{c}).$ {Consequently, when $\xi_c=\mathrm{O}(1),$ even for $j=j^*$ in (\ref{j_assumption}),  (\ref{eq_rateassumption}) only requires that $c(n^{-1/2+1/(2(\tau-1))+n^{-1+1/(2(\tau-1))+2/(\tau+1)}})\ll 1.$ In particular, if we set $c=\mathrm{O}(n^{a})$ for some small constant $0<a<1/2,$ we only need $\tau>1+\frac{1}{1-2a}.$  

Second, for the rate $\Psi(j,c)$ in (\ref{eq_consistentrate}) and (\ref{eq_psiparameter}), when $\xi_c=\mathrm{O}(1)$ and $\zeta_c=\mathrm{O}(\sqrt{c}),$ even for $j=j^*,$  it reads $c^{-d}(1+cn^{-1/2+1/(\tau-1)})+cn^{-1/2+1/(2(\tau-1))}+cn^{-3/2+3/(2(\tau-1))}$. Therefore, for sufficiently large $d$ and $\tau,$ it has an order of $c n^{-1/2+\epsilon},$ for some small constant $\epsilon>0.$}

\end{rem}
\subsection{Asymptotic normality and power analysis for the proposed statistics}\label{sec_inferencepacf}
In this section, we study the accuracy and power of the proposed statistics $T_1$ in (\ref{eq_teststatisticone}) and $T_2$ in (\ref{eq_defnt2}). We first prepare some notations. Recall (\ref{eq_jorderlocal}). Following the conventions below (\ref{eq_residualestimate}), for $\ell=j$ regarding $T_1$ and $\ell=\mathsf{h}$ regarding $T_2,$ we denote \vspace*{-0.1cm}
\begin{equation}\label{eq_wiform}
\bm{w}_i=\bm{x}_{\ell, i} \epsilon_{i,\ell}, \ i>\ell,
\end{equation}
 where we recall $\bm{x}_{\ell,i}=(x_{i-1}, \cdots, x_{i-\ell})^* \in \mathbb{R}^{\ell}.$ According to (3.17) of \cite{DZar} or Lemma  3.1 of \cite{ding2021simultaneous}, we see that $\bm{w}_i$ has a physical representation in the sense that for some measurable function $\mathbf{V},$ we have \vspace*{-0.1cm}
\begin{equation}\label{eq_defnwi}
\bm{w}_i=\mathbf{V}(\frac{i}{n}, \mathcal{F}_i), \ i>\ell,
\end{equation}
{where $\mathcal{F}_i=(\cdots, \eta_{i-1}, \eta_i)$ and $\eta_i, \ i  \in \mathbb{Z}$ are i.i.d centered random variables.} Denote the long-run covariance matrix of $\bm{w}_i$ as 
$
\Pi(t)=\sum_{j=-\infty}^{\infty} \text{Cov} \Big( \mathbf{V}(t, \mathcal{F}_0), \mathbf{V}(t, \mathcal{F}_j) \Big),
$
and the integrated long-run covariance matrix as 
$
\Pi= \int_{0}^1 \Pi(t) \otimes (\mathbf{B}(t) \mathbf{B}^*(t)) \mathrm{d}t. 
$  In what follows, to ease our discussion, we assume that $c$ is of the form \vspace*{-0.1cm}
\begin{equation}\label{eq_cform}
c=\rO(n^{a}), \ 0<a<1. 
\end{equation}
Armed with the above notations, we now proceed to provide the theoretical properties of the statistic $T_1.$ Recall the matrix $\mathbf{M}$ in (\ref{eq_finalfinalsamplesamplestatistic}). For $s \in \mathbb{N},$ define \vspace*{-0.1cm}
\begin{equation}\label{eq_defingk}
g_{1,s}=(\text{Tr}[\Pi^{1/2} \overline{\Sigma}_{01}^{-1} \mathbf{M} \overline{\Sigma}_{01}^{-1} \Pi^{1/2}]^s)^{1/s}, \ \overline{\Sigma}_{01}=\begin{pmatrix}
\mathbf{I}_c & \bm{0} \\
\bm{0} & \Sigma_{01}
\end{pmatrix},
\end{equation}
where $\Sigma_{01}:=\int_0^1 \Sigma^{j}(t) \otimes (\mathbf{B}(t) \mathbf{B}^*(t))\mathrm{d}t$ and $\Pi$ is defined from (\ref{eq_defnwi}) using $\ell=j.$ \vspace*{4pt}

\begin{thm}\label{lem_controltbtrho0} Suppose the assumptions of Theorem \ref{thm_consistency} hold.  Then we have that
\begin{enumerate}
\item Suppose Assumption \ref{assu_ARbapproximation} of our supplement holds and \vspace*{-0.1cm}
\begin{equation}\label{eq_anotherbound1}
\Psi(j,c)=\mathrm{o}(1), 
\end{equation}
we have that when $\mathbf{H}_{01}$ in (\ref{eq_firstnull}) holds \vspace*{-0.1cm}
\begin{equation}\label{eq_distribution11}
\mathbb{T}_1:=\frac{n T_{1}-g_{1,1}}{g_{1,2}} \Rightarrow \mathcal{N}(0,2). 
\end{equation}
\item  When $\mathbf{H}_{a1}$ in (\ref{eq_firstnull}) holds in the sense that \vspace*{-0.1cm} 
\begin{equation}\label{eq_alternativespecificform}
\int_0^1 \rho_j(t)^2 \mathrm{d} t >\mathfrak{C}, \ \text{where} \ \mathfrak{C}:=C_{\alpha} \frac{\sqrt{j c}}{n},
\end{equation}
where $C_\alpha \equiv C_\alpha(n) \rightarrow  \infty$ as $n \rightarrow \infty,$ assuming (\ref{eq_anotherbound}),  then we have that for any $\alpha \in (0,1)$ \vspace*{-0.1cm}
\begin{equation*}
\mathbb{P} \left( \left| \mathbb{T}_1 \right| \geq \sqrt{2} \mathcal{Z}_{1-\alpha} \right) \rightarrow 1, \ n \rightarrow \infty,
\end{equation*}
where $\mathcal{Z}_{1-\alpha}$ is the $(1-\alpha)$th quantile of the standard Gaussian distribution.
\end{enumerate}
 \end{thm}

The above theorem establishes the asymptotic normality for our proposed statistic $T_1$ in (\ref{eq_teststatisticone}) concerning (\ref{eq_firstnull}).  Moreover,  (\ref{eq_alternativespecificform}) shows that our proposed statistic can have asymptotic power one under weak local alternative.  Then we provide the theoretical properties of the statistic $T_2.$  For $s \in \mathbb{N},$ define \vspace*{-0.1cm}
\begin{equation}\label{eq_defingkaa}
g_{2,s}=(\text{Tr}[\Pi^{1/2} \overline{\Sigma}_{02}^{-2}   \Pi^{1/2}]^s)^{1/s}, \ \overline{\Sigma}_{02}=\begin{pmatrix}
\mathbf{I}_c & \bm{0} \\
\bm{0} & \Sigma_{02}
\end{pmatrix},
\end{equation}
where $\Sigma_{02}:=\int_0^1 \Sigma^{\mathsf{h}}(t) \otimes (\mathbf{B}(t) \mathbf{B}^*(t))\mathrm{d}t$ and $\Pi$ is defined from (\ref{eq_defnwi}) using $\ell=\mathsf{h}.$ 
\vspace*{4pt}

\begin{thm}\label{lem_controltbtrho} 
Suppose the assumptions of Theorem \ref{thm_consistency} hold.  Then we have that
\begin{enumerate}
\item Suppose Assumption \ref{assu_ARbapproximation} of our supplement holds and \vspace*{-0.1cm}
\begin{equation}\label{eq_anotherbound}
\mathsf{h}\left(\frac{\mathsf{h}^4}{n}+\mathsf{h}^2 \Psi^2(\mathsf{h},c) \right)=\mathrm{o}(1), 
\end{equation}
when $\mathbf{H}_{02}$ in (\ref{eq_testnotwo}) holds, we have that \vspace*{-0.1cm}
\begin{equation}\label{eq_distribution11}
\mathbb{T}_2:=\frac{n T_{2}-g_{2,1}}{g_{2,2}} \Rightarrow \mathcal{N}(0,2). 
\end{equation}
\item  Recall (\ref{eq_Tbp}). Suppose (\ref{eq_anotherbound}) holds. Then when $\mathbf{H}_{02}$ in (\ref{eq_testnotwo}) holds, we have that \vspace*{-0.1cm}
\begin{equation}\label{eq_closenessstatistics}
T_{\text{BP}}=T_{2}+\mathrm{O}_{\mathbb{P}} \left( \mathsf{h} \Psi^2(\mathsf{h},c)+ \frac{\mathsf{h}^5}{n^2}+\sum_{k=1}^{h_2} \Psi^2(k,c)\right).
\end{equation}
Consequently, if we further assume that \vspace*{-0.1cm}
\begin{equation}\label{eq_anotherassumption}
\frac{n\left[ \mathsf{h} \Psi^2(\mathsf{h},c)+ \frac{\mathsf{h}^5}{n^2}+\sum_{k=1}^{\mathsf{h}} \Psi^2(k,c) \right]}{\sqrt{\mathsf{h}c}}=\mathrm{o}(1),  
\end{equation}
then \vspace*{-0.1cm}
\begin{equation}\label{eq_distribution}
\frac{n T_{\text{BP}}-g_{2,1}}{g_{2,2}} \Rightarrow \mathcal{N}(0,2). 
\end{equation}
\item When $\mathbf{H}_{a2}$ in (\ref{eq_testnotwo}) holds in the sense that \vspace*{-0.1cm}
\begin{equation}\label{eq_alternativespecificform2ndapplication}
\sum_{k=1}^{\mathsf{h}} \int_0^1 \rho_k^2(t) \mathrm{d} t >\mathfrak{C}, \ \text{where} \ \mathfrak{C}:=C_{\alpha}  \frac{\sqrt{\mathsf{h} c}}{n},
\end{equation}
where $C_\alpha \equiv C_\alpha(n) \rightarrow  \infty$ as $n \rightarrow \infty,$ assuming the assumptions of parts 1  and 2 hold, then we have that for any $\alpha \in (0,1)$ \vspace*{-0.1cm}
\begin{equation*}
\mathbb{P} \left( \left|\mathbb{T}_2 \right| \geq \sqrt{2} \mathcal{Z}_{1-\alpha} \right) \rightarrow 1, \ n \rightarrow \infty,
\end{equation*}
where $\mathcal{Z}_{1-\alpha}$ is the $(1-\alpha)$th quantile of the standard Gaussian distribution.
\end{enumerate}
\end{thm}

 Theorem \ref{lem_controltbtrho} implies that $T_{\text{BP}}$ and $T_2$ have the same Gaussian distribution and power performance asymptotically so that we can directly use $T_{2}$ which only requires one OLS regression for estimation. {We point out that (\ref{eq_anotherbound}) and (\ref{eq_anotherassumption}) basically impose some upper bound conditions for $\mathsf{h}.$ As discussed in Remark \ref{rem_consist}, if one chooses Fourier or orthogonal polynomial as the basis functions, when $d$ and $\tau$ are large enough, we require $\mathsf{h}\ll n^{1/5}$ to guarantee (\ref{eq_anotherbound}). Analogously, (\ref{eq_anotherassumption}) requires that $\mathsf{h}\ll \min\{c^{1/5}, n^{1/5}\}.$ Recall \eqref{j_assumption}. In particular, $j^*=o(n^{1/5})$ when $\tau> 2.5$. Note that $\mathsf{h}$ is only required to be larger or equal to $j^*$. As a result, the above constraints on $\mathsf{h}$ is not restrictive.}
\begin{rem}\label{rem_keykeykey}
We add a remark on the assumptions of the parameters. First, (\ref{eq_anotherbound1}) and (\ref{eq_anotherbound}) are mainly used to guarantee that $nT_k, k=1,2,$ can be written as a quadratic form in terms of $\bm{w}_i$ in (\ref{eq_wiform}) under the null hypotheses. More specifically, denote $\mathbf{X}=\frac{1}{\sqrt{n}} \sum_{i=\ell+1}^n (\bm{w}_i \otimes \mathbf{B}(i/n)),$ we can conclude from \cite{DZar} that
\begin{equation}\label{eq_reducedform}
n T_1=\mathbf{X}^* \overline{\Sigma}_{01}^{-1} \mathbf{M} \overline{\Sigma}_{01}^{-1} \mathbf{X}+\mathrm{o}_{\mathbb{P}}(1), \  \  n T_2=\mathbf{X}^* \overline{\Sigma}_{02}^{-2} \mathbf{X}+\mathrm{o}_{\mathbb{P}}(1).
\end{equation}
Second, (\ref{eq_closenessstatistics}) guarantees the difference between $T_{\text{BP}}$ and $T_2$ is negligible when $\mathbf{H}_{02}$ in (\ref{eq_testnotwo}) holds. 
\end{rem}
 Before concluding this section, we summarize the properties of $\widehat{\mathcal{T}}_{k}, k=1,2,$ in (\ref{eq_finalfinalsamplesamplestatistic}) which explains the validity for the bootstrap procedure. The motivation comes from the arguments as in Remark \ref{rem_keykeykey} that the statistic $n T_k, k=1,2,$ are essentially quadratic forms of the locally stationary vector (c.f. (\ref{eq_wiform})) and the statistics are asymptotically Gaussian. Consequently, it is possible to mimic their asymptotic distribution using a multiplier bootstrap procedure. 
\vspace*{4pt}

\begin{cor}\label{cor_boostrap} Suppose the assumptions of Theorems \ref{lem_controltbtrho0} and \ref{lem_controltbtrho} hold. Moreover, assume the time series has finite fourth moment and $m \rightarrow \infty$ as $n \rightarrow \infty$ and \vspace*{-0.1cm}
\begin{equation*}
\sqrt{\mathsf{h}} \zeta_c^2 c^{-1/2} \left(\frac{m}{n}+\frac{1}{m} \right)=\mathrm{o}(1).
\end{equation*}
Then {there exists a sequence of sets $\mathcal{A}_n$ such that $\mathbb{P}(\mathcal{A}_n)=1-\mathrm{o}(1)$ and under the event $\mathcal{A}_n,$ we have that conditional on the data $\{x_i\},$} the results in Theorems \ref{lem_controltbtrho0} and \ref{lem_controltbtrho} still hold by replacing $nT_{k}$ with $\widehat{\mathcal{T}}_k, k=1,2.$  
\end{cor}

\vspace*{-0.5cm}
\section{Numerical simulations and real data analysis}\label{sec_simulation}

In this section, we use some numerical simulations and a real data analysis to illustrate the usefulness of our estimation and inference procedures of the PACFs.  All the calculations, implementations and plots can be done using a few lines of coding with our $\mathtt{R}$ package $\mathtt{Sie2nts}.$
\subsection{Numerical simulations}

In this section, we conduct numerical simulations to illustrate the usefulness of our methodologies using both stationary and non-stationary models. {Due to space constraint, we focus on reporting the results of AR type models in the following. Additional simulation results on other types of models can be found in Section \ref{sec_additionalsimulationresults} of our supplement.}

In what follows, for some $\delta_1, \delta_2, \in [0, 0.5],$ we consider the stationary AR(2) process
\begin{equation}\label{eq_stationarymodel}
x_i=\delta_1 x_{i-1}+\delta_2 x_{i-2}+ \epsilon_i,  
\end{equation} 
and the locally stationary AR(2) process 
\begin{equation}\label{eq_nonstationarymodel}
x_i=\delta_1 \sin (2 \pi i/n)x_{i-1}+\delta_2 \cos(2 \pi i/n) x_{i-2}+\left( 0.4+0.4|\sin (2 \pi i/n)| \right)\epsilon_i,
\end{equation}
where $\epsilon_i, 1 \leq i \leq n,$ are i.i.d. standard Gaussian random variables. Note that when $\delta_1=\delta_2=0,$ (\ref{eq_stationarymodel}) reduces to a standard white noise process and (\ref{eq_nonstationarymodel}) reduces to a time-varying white noise process.  
\vspace*{-0.3cm}
\subsubsection{Estimation of PACFs}
We first estimate the PACFs using the sieve method as introduced in Section \ref{sec_estimation}. For concreteness and due to space constraint, regarding  $\delta_1$ and $\delta_2$ in (\ref{eq_stationarymodel}) and (\ref{eq_nonstationarymodel}), we only report the results for the choice $\delta_1=0.5$ and $\delta_2=0.3.$ Note that similar results and conclusions can also be made for other choices.

\begin{figure}[!ht]
\hspace*{0cm}
\begin{subfigure}{0.55\textwidth}
\includegraphics[width=6.5cm,height=6cm]{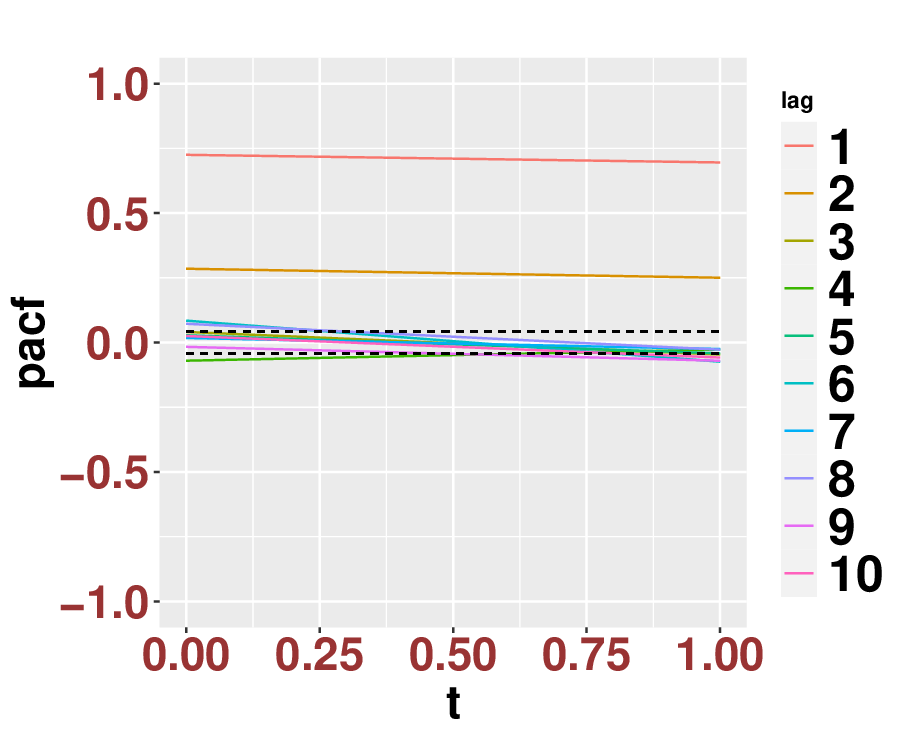}
\caption{PACFs for model (\ref{eq_stationarymodel}).}\label{subfig_nullorthogonaltypei200}
\end{subfigure}
\hspace{0cm}
\begin{subfigure}{0.55\textwidth}
\includegraphics[width=6.5cm,height=6cm]{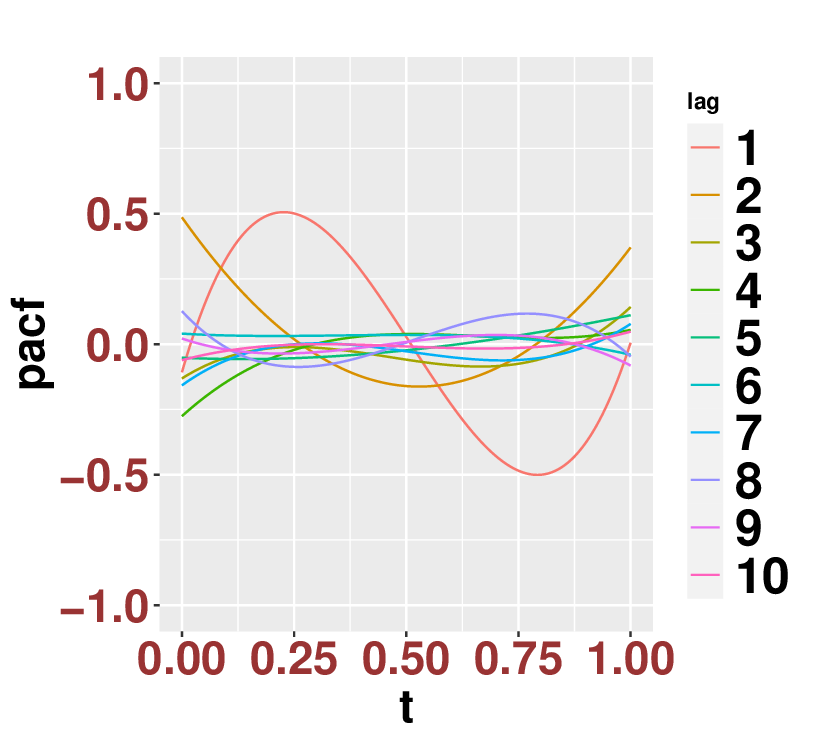}
\caption{PACFs for model (\ref{eq_nonstationarymodel}).}\label{subfig_nullorthogonaltypei500}
\end{subfigure}
\vspace*{-0.8cm}
\caption{{ {Typical sample PACF} plots (i.e., $\widehat{\rho}_j(t)$ in (\ref{eq_pacfestimate})) for models (\ref{eq_stationarymodel}) and (\ref{eq_nonstationarymodel}). Here $n=600$ and the plots  can be generated using the function $\mathtt{sie.auto.plot}$ from our $\mathtt{R}$ package $\mathtt{Sie2nts}$. }}
\label{fig_pacfplot}
\end{figure}

In Figure \ref{fig_pacfplot}, we estimate the PACFs of the first 10 lags for both models (\ref{eq_stationarymodel}) and (\ref{eq_nonstationarymodel}) using the estimators in (\ref{eq_pacfestimate}). We use the Legendre polynomials as the basis functions and the number $c$ can be chosen using the cross validation method as described in Section \ref{sec:choiceparameter}.  The computations of the PACFs can be obtained directly using the function $\mathtt{sie.plot.pacf}$ from our $\mathtt{R}$ package. 
From these plots, we can see that our method applies to both models and obtain reasonably accurate estimates. According to the cut-off properties of the PACFs of the AR models,  these plots have also implied that the time series may be generated from some AR(2) models. {In Section \ref{sec_simulationcomparison} of our supplement, we compare our proposed method with the ones introduced in \cite{10.1214/20-EJS1748} and find that our method is generally more accurate than \cite{10.1214/20-EJS1748} in terms of mean integrated squared error. The main reason is that the sieve method is adaptive to the smoothness of the covariance structure and has less boundary issues compared with the kernel method. In addition, more simulation results on other types of models can be found in Section \ref{sec_simulationmore} of our supplement and similar conclusions can also be made. } 


\vspace*{-0.3cm}
\subsubsection{Inference of PACFs}
In this section, we examine the accuracy and sensitivity of our proposed Algorithm \ref{alg:boostrapping} when applied to testing (\ref{eq_firstnull}) and (\ref{eq_testnotwo}). We first investigate the accuracy. To test (\ref{eq_firstnull}) for some individual PACF at lag $j$, in the context of (\ref{eq_stationarymodel}) and (\ref{eq_nonstationarymodel}), we consider the following four settings: (1).  $\delta_1=0.5, \delta_2=0,$ and $j=2$; (2). $\delta_1=0.5, \delta_2=0,$ and $j=4$; (3). $\delta_1=\delta_2=0.3,$ and $j=3;$ (4). $\delta_1=\delta_2=0.3,$ and $j=5.$ Moreover, to test (\ref{eq_testnotwo}) for the white noise, we consider the following setting in the context of (\ref{eq_stationarymodel}) and (\ref{eq_nonstationarymodel}): (5).  $\delta_1=\delta_2=0.$

In Table \ref{table_typei}, we report the simulated type I error rates for all the above null settings for three different types of basis functions when $n=600.$  We can see that our Algorithm \ref{alg:boostrapping} is quite accurate for both tests (\ref{eq_firstnull}) and (\ref{eq_testnotwo}).

\begin{table}[!ht]
\begin{center}
\setlength\arrayrulewidth{1pt}
\renewcommand{\arraystretch}{1.3}
{\fontsize{9}{9}\selectfont 
\begin{tabular}{|c|ccccc|ccccc|}
\hline
      & \multicolumn{5}{c|}{$\alpha=0.1$}                                                                                                                       & \multicolumn{5}{c|}{$\alpha=0.05$}                                                                                                                        \\ \hline
Basis/Setting & \multicolumn{1}{c|}{(1)} & \multicolumn{1}{c|}{(2)} & \multicolumn{1}{c|}{(3)} & \multicolumn{1}{c|}{(4)} & \multicolumn{1}{c|}{(5)}   & \multicolumn{1}{c|}{(1)} & \multicolumn{1}{c|}{(2)} & \multicolumn{1}{c|}{(3)} & \multicolumn{1}{c|}{(4)} & \multicolumn{1}{c|}{(5)}  \\ 
\hline
     & \multicolumn{10}{c|}{Model (\ref{eq_stationarymodel})}                                                                                                                                                                                                                                                                                          \\
   \hline
Fourier     &          0.108  & 0.096                          &                          0.108 &         0.11                 &                          0.098 &                    0.048    & 0.059   &                                     0.061 &       0.054 & 0.049  \\
Legendre    &     0.109       & 0.1                          &                          0.096 &   0.103                        &           0.108                   &       0.06    & 0.053         &  0.055  &                                     0.064 &                            0.048  \\
Daubechies-9    &  0.096 & 0.102   & 0.11                        &      0.101             &         0.095              & 0.058                 &0.061 &  0.054  &          0.058                            &                 0.043    \\
\hline
      & \multicolumn{10}{c|}{Model (\ref{eq_nonstationarymodel})}                                                                                                                                                                                                                                                                                         \\
       \hline
Fourier     &            0.098      & 0.11                 &                          0.108 &      0.099                     &  0.107   &                         0.048 & 0.062 & 0.057   &    0.048                                  &              0.039                 \\
Legendre     &    0.103    & 0.094    &                  0.103                                &          0.095                &                       0.113  & 0.055  & 0.053                       & 0.061    &        0.048                              &                            0.053  \\
Daubechies-9     &  0.093  & 0.09    &      0.093                   &                         0.104   & 0.1  &          0.052               &  0.058 & 0.058 &                                 0.047  &                  0.046  \\
 \hline
\end{tabular}
}
\end{center}
\vspace*{-0.4cm}
\caption{ { Simulated type I error rates.  The results are reported based on 1,000 simulations. We can see that our multiplier bootstrap procedure is reasonably accurate for both $\alpha=0.1$ and $\alpha=0.05$. In our $\mathtt{R}$ package $\mathtt{Sie2nts},$ we collect 31 different types basis functions. 
For visualization of various types of basis functions, we can use the $\mathtt{bs.gene}$ and $\mathtt{bs.plot}$ from our $\mathtt{R}$ package $\mathtt{Sie2nts}.$ The computations of the $p$-values can be obtained directly using the function $\mathtt{auto.pacf.test}$ from the package.} }
\label{table_typei}
\end{table}

Then we study the power. For definiteness, we report the results for the white noise test that  {
\begin{equation}\label{eq_nullandalternativerespectively}
\mathbf{H}_0: \ \{x_i\} \ \text{is a white noise process} \ \ \text{Vs} \ \ \mathbf{H}_a: \ \{x_i\} \ \text{is not a white noise process.}
\end{equation} }
More specifically, in view of the concerned models (\ref{eq_stationarymodel}) and (\ref{eq_nonstationarymodel}), {$\mathbf{H}_0$ in (\ref{eq_nullandalternativerespectively}) is equivalent to $\delta_1=\delta_2=0$ while $\mathbf{H}_a$ is an AR(1) alternative that $\delta_2=0, \delta_1>0.$} In Figure \ref{fig_powerplot} below, we study the power of our proposed method when $\delta_1$ deviates away from zero. We can conclude that our proposed tests are reasonably powerful once the alternative deviates from the null. 

\begin{figure}[!ht]
\hspace*{0cm}
\begin{subfigure}{0.55\textwidth}
\includegraphics[width=6.5cm,height=6cm]{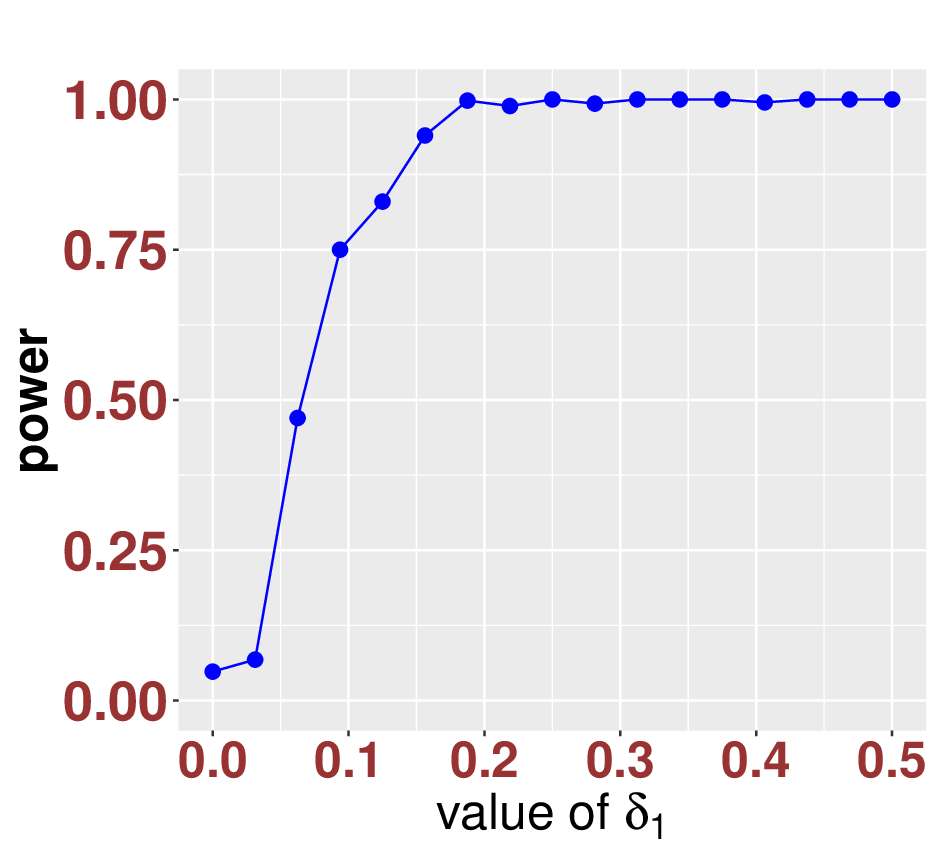}
\caption{Power for model (\ref{eq_stationarymodel}).}\label{subfig_nullorthogonaltypei200}
\end{subfigure}
\hspace{0cm}
\begin{subfigure}{0.55\textwidth}
\includegraphics[width=6.5cm,height=6cm]{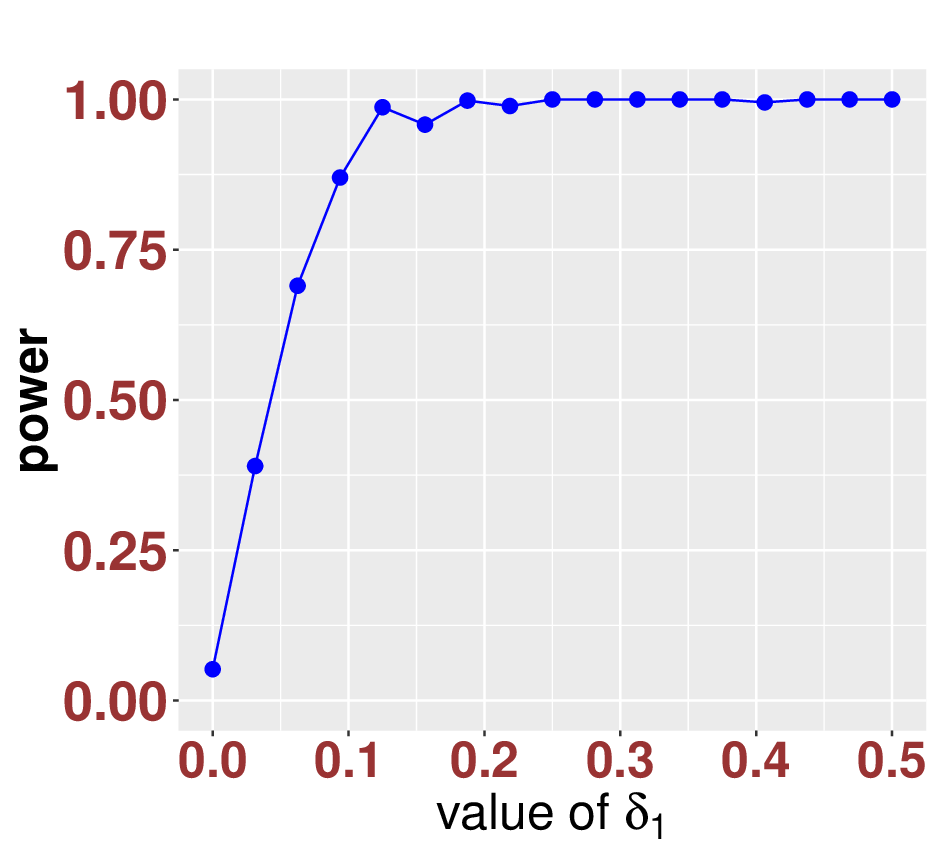}
\caption{Power for model (\ref{eq_nonstationarymodel}).}\label{subfig_nullorthogonaltypei500}
\end{subfigure}
\vspace*{-0.8cm}
\caption{{ Power for models (\ref{eq_stationarymodel}) and (\ref{eq_nonstationarymodel}) under the alternative of (\ref{eq_nullandalternativerespectively}). Here the type I error rate $\alpha=0.05$ and $n=600$. We use the Legendre polynomials as the basis functions and the computations of the $p$-values can be obtained directly using the function $\mathtt{auto.pacf.test}$ from our $\mathtt{R}$ package $\mathtt{Sie2nts}$. The results are reported based 1,000 repetitions.}}
\label{fig_powerplot}
\end{figure}

Finally, for further visualization, using  (\ref{eq_nonstationarymodel}) as an example, in Figure \ref{fig_pvalueplot} below, we provide two typical plots of the $p$-values associated with each lag $j$ with (\ref{eq_firstnull}) under both the null and alternative as in (\ref{eq_nullandalternativerespectively}). Here for the alternative we set $\delta_1=0.5$ and $\delta_2=0$ in (\ref{eq_nonstationarymodel}). From the plots we can easily distinguish the null and alternative. Moreover, it is clear that the plots suggest that the null is a white noise process and the alternative is an AR(1) process. {Additionally, more simulation results on other types of models can be found in Section \ref{sec_simulationmore} of our supplement and similar conclusions can also be made.} 

\begin{figure}[!ht]
\hspace*{0cm}
\begin{subfigure}{0.55\textwidth}
\includegraphics[width=6.5cm,height=6cm]{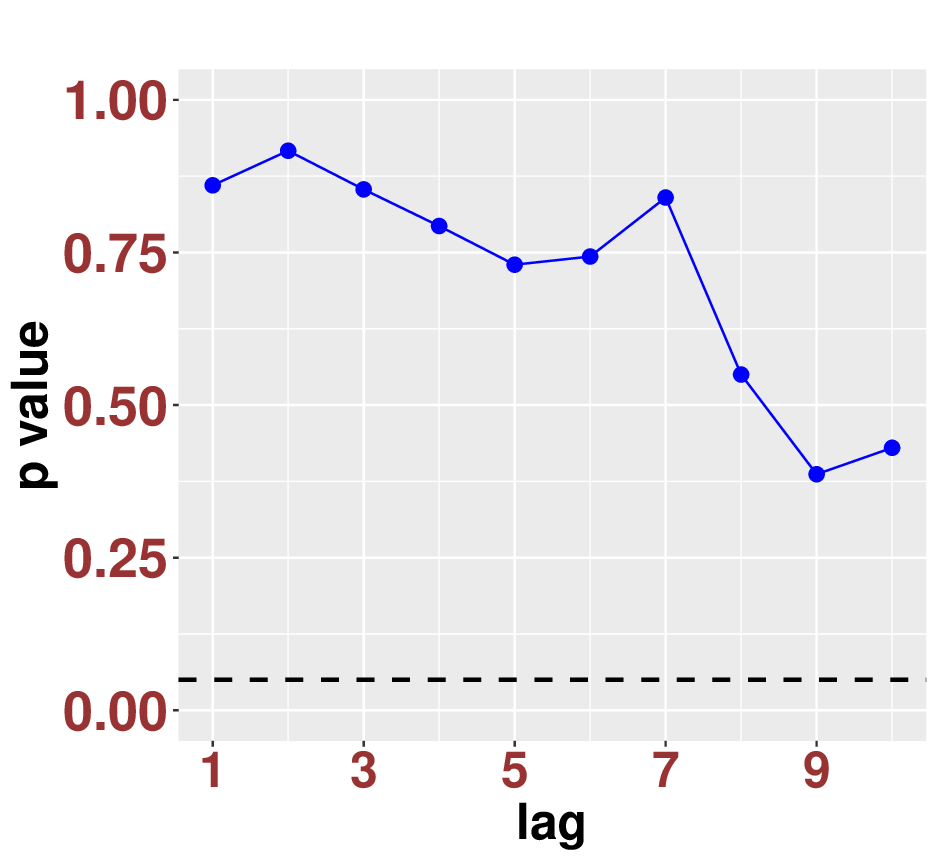}
\caption{$p$-values under $\mathbf{H}_0$ in (\ref{eq_nullandalternativerespectively}).}\label{subfig_nullorthogonaltypei200}
\end{subfigure}
\hspace{0cm}
\begin{subfigure}{0.55\textwidth}
\includegraphics[width=6.5cm,height=6cm]{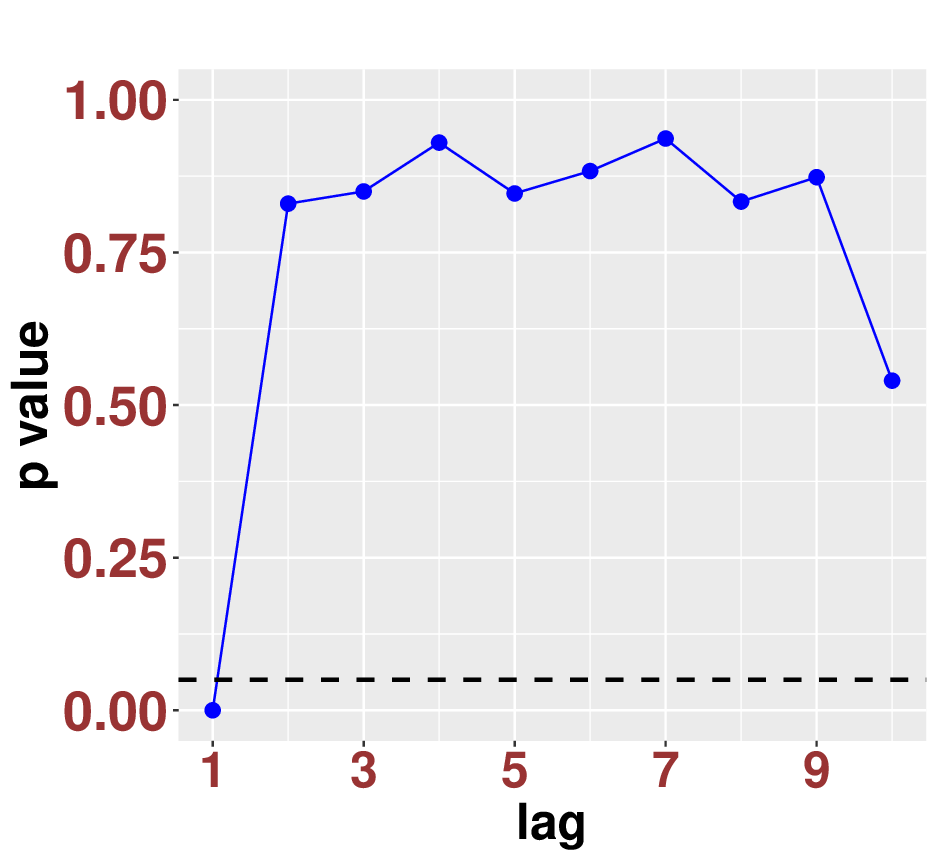}
\caption{$p$-values under $\mathbf{H}_a$ in (\ref{eq_nullandalternativerespectively}).}\label{subfig_nullorthogonaltypei500}
\end{subfigure}
\vspace*{-0.8cm}
\caption{{ Typical $p$-value plots for different lags for model (\ref{eq_nonstationarymodel}) under null and alternative of (\ref{eq_nullandalternativerespectively}). Here the type one error $\alpha=0.05,$ $n=600$, $\delta_1=0.5$ for the alternative. We use the Legendre polynomials as the basis functions and the computations of the $p$-values can be obtained directly using the function $\mathtt{auto.pacf.test}$ from our $\mathtt{R}$ package $\mathtt{Sie2nts}$.  }}
\label{fig_pvalueplot}
\end{figure}


\vspace*{-0.3cm}

\subsection{Real data analysis}

In this section, we apply our methods to the monthly Euro-Dollar exchange rate data set which has also been considered in \cite{10.1214/20-EJS1748}. The data set can be downloaded from EuroStat at \url{https://ec.europa.eu/eurostat/web/products-datasets/-/ei_mfrt_m}. We analyze the exchange rates from January 1999 until October 2017.

\begin{figure}[htp]
\begin{subfigure}{0.55\textwidth}
\includegraphics[width=6.5cm,height=6.1cm]{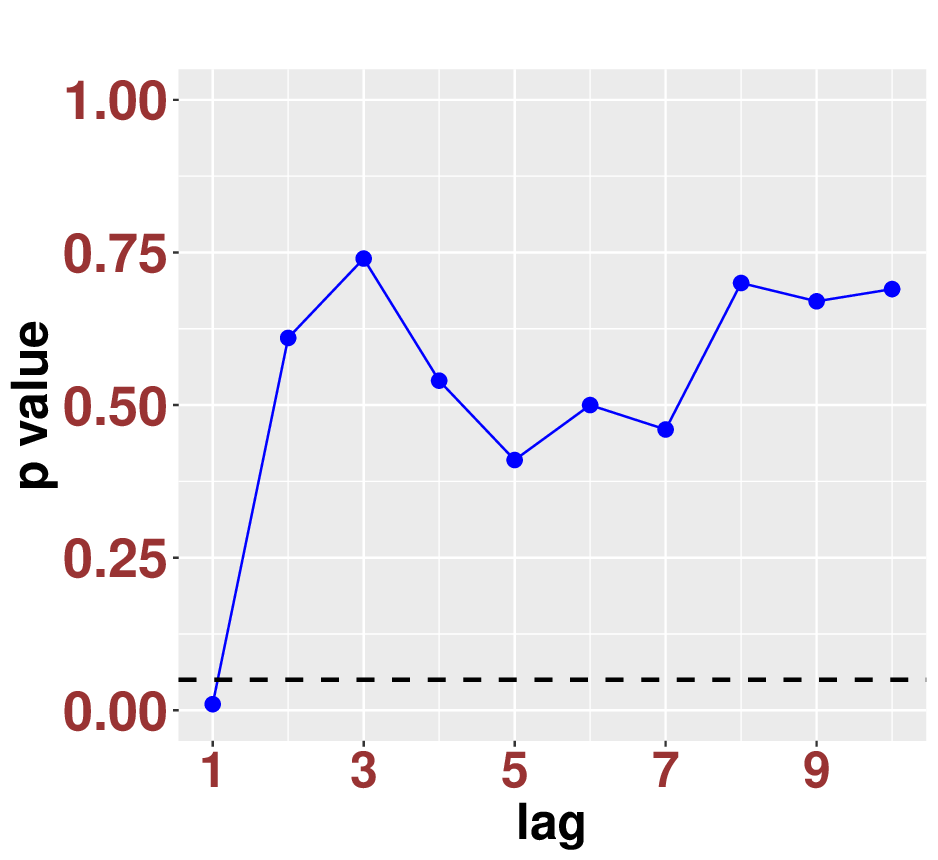}
\end{subfigure}
\begin{subfigure}{0.55\textwidth}
\includegraphics[width=6.5cm,height=6cm]{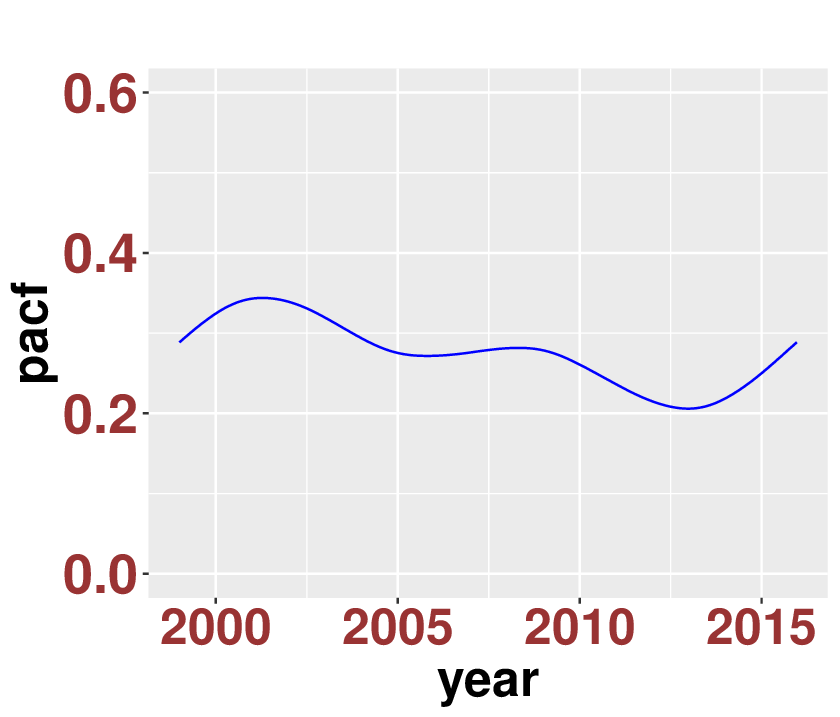}
\end{subfigure}
\vspace*{-0.8cm}
\caption{Euro-Dollar exchange rate data. Here we use the Daubechies-9 wavelet as the basis functions. Left panel records the $p$-values associated with different lags and right panel is the estimation of the PACF at the first lag.}
\label{fig_data2}
\end{figure} 

Following the traditions of financial data analysis, we consider the log returns of the exchange rate. Moreover, we apply our methods to estimate the PACFs and conduct inference on them. As can be seen from the  left panel of Figure \ref{fig_data2}, based on our inference, a white-noise-driven AR(1) model will be useful to model the data. This agrees with the findings in \cite{10.1214/20-EJS1748}. {In fact, we can conduct the white noise test for the original time series as in (\ref{eq_testnotwo}) using Algorithm \ref{alg:boostrapping} and find the $p$-value is $0.013.$ This shows that the time series is not a white noise. Furthermore, after fitting the time series with a time-varying  AR(1) model (note that for AR(1) model, the coefficient  is identical to $\widehat{\rho}_1(t)$), we can conduct the white noise test for the residuals and then find the $p$-value is $0.54$. This confirms that the AR(1) model is likely to be appropriate.} 

In the right panel of Figure \ref{fig_data2}, we  provide the estimation of the PACF of lag one which is also the coefficient of a time-varying AR(1) model. It can be seen that the even though the PACF is relatively stable, it experiences some smooth changes over time. {In fact, one can directly generalize the significance test $T_1$ into a test of constancy by considering
$T_1^*:=\int_0^1 [\widehat{\rho}_1(t)-\int_0^1\widehat{\rho}_1(s) \mathrm{d} s]^2 \mathrm{d} t. $
The multiplier bootstrap of $T_1^*$ yields a $p$-value of $0.031$, indicating that $\rho_1(t)$ is likely to be time-varying.}
This concludes that a locally stationary AR(1) model may be more useful for the model fitting.   {Due to the pronounced responsiveness of Euro-Dollar exchange rates to the global economy, our plots serve as reasonably accurate reflections of prevailing global economic conditions. For instance, leading up to the global financial crisis (2005-2007), the PACF displays a consistent pattern. Subsequently, from 2008 onwards, there is a gradual decline observed, hitting its lowest point around 2013, coinciding with the recognized period of the financial crisis. Following this, the PACF demonstrates a resurgence. Such a visual representation offers valuable insights into the evolving dynamics of Euro-Dollar exchange rates, aiding in a deeper comprehension of their temporal behavior.}


\vspace*{-0.3cm}

\section*{Supplementary file}
In the supplement file, we provide the technical proofs, some additional remarks, practical methods for choosing the tuning parameters and additional simulation results. 

 \setcounter{section}{0}
 \setcounter{table}{0}
        \setcounter{figure}{0}
\numberwithin{equation}{section}
\renewcommand{\thesection}{\Alph{section}}
\renewcommand{\thetable}{\Alph{table}}
\renewcommand{\thefigure}{\Alph{figure}}  
  
\section{Technical proofs}\label{sec_techinicalproof}

\subsection{Proofs of Section \ref{sec_characterization}}

In this subsection, we prove Theorem \ref{lem_pacf}. The strategies and ideas are similar to those of Theorems 2.4 and 2.11 of \cite{DZar}. We focus on explaining the main differences. 
\vspace*{4pt}

\begin{Proof}[\bf Proof of Theorem \ref{lem_pacf}]
For the first part of the proof, it is analogous to the proof of equation (2.6) of \cite{DZar}. We only sketch the proof. Due to similarity, we only prove the first control in (\ref{eq_decayrate}). Recall (\ref{eq_jorderlocal}). Set $\bm{\phi}_{j}(i)=(\phi_{ij,1},\cdots, \phi_{ij,j})^* \in \mathbb{R}^j.$ It is easy to see from (1) of Assumption \ref{assum_mainassumption} and Yule-Walker's equation that 
\begin{equation}\label{eq_definitionphji}
\bm{\phi}_j(i)=\Gamma_j^{-1}(i) \bm{\nu}_j(i),
\end{equation}  
where we denoted that 
\begin{equation*}
\Gamma_j(i)=\operatorname{Cov} \left( \bm{x}_j(i), \bm{x}_j(i) \right) \in \mathbb{R}^{j \times j}, \  \bm{\nu}_j(i)=\operatorname{Cov}\left( \bm{x}_j(i), x_i \right) \in \mathbb{R}^{j},
\end{equation*}
where $\bm{x}_j(i)=(x_{i-1}, \cdots, x_{i-j}) \in \mathbb{R}^j.$  For the rest of the proof, we can follow lines of that of Theorem 2.4 of \cite{DZar} verbatim. More specifically, according to (\ref{eq_definitionphji}), we find that in order to study $\rho_{i,j} \equiv \phi_{ij,j},$ it suffices to control the entries of the $j$th row of  $\Gamma_j^{-1}(i)$ and all the entries of $\bm{\nu}_j(i).$ For $\Gamma_j^{-1}(i),$ we denote the $j \times j$ symmetric banded matrix $\Gamma_j^s(i)$ as 
\begin{equation*}
(\Gamma_j^s(i))_{kl}=
\begin{cases} 
(\Gamma_j(i))_{kl}, & |k-l| \leq \frac{j}{K \log j}; \\
0, & \text{otherwise}.
\end{cases}
\end{equation*}
Here $K>0$ is some large constant. Using (1) and (2) of Assumption \ref{assum_mainassumption}, by a discussion similar to equation (D.7) of \cite{DZar}, we find that for some constant $C>0$
\begin{equation*}
|\rho_{i,j}-\phi_{ij,j}^s| \leq C j^{1-\tau} (K \log j)^{\tau-1},
\end{equation*}
where $\bm{\phi}_j(i)^s=(\phi_{ij,1}^s, \cdots, \phi_{ij,j}^s)^* \in \mathbb{R}^j$ is defined via $\bm{\phi}_j^s(i)=(\Gamma_j(i)^s)^{-1} \bm{\nu}_j(i).$ Moreover, using (1) of Assumption \ref{assum_mainassumption}, according to the discussions below equation (D.8) of \cite{DZar}, we find that for some constant $C_1>0,$
\begin{equation*}
|\phi_{ij,j}^s| \leq C_1(j/(\log j+1))^{-\tau+1}. 
\end{equation*}
Combining the above controls, we complete the proof of (\ref{eq_decayrate}).

For the second part of the proof, for the smoothness that $\rho_j(t) \in C^d([0,1]) $ follows from Lemma 3.1 of \cite{DZ1}. To see (\ref{eq_smoothapproximation}), by (\ref{eq_definitionphji}) and (\ref{eq_defnphjt}), using Cauchy-Schwarz inequality, we have that 
\begin{align}\label{eq_normcontrol}
\left | \bm{\phi}_j(i)-\bm{\phi}_j(i/n) \right |  \leq \left\| \Gamma_j^{-1}(i) \right\| \left| \bm{\nu}_j(i)-\bm{\nu}_j(i/n) \right|+\left\|  \Gamma_j^{-1}(i)-\Gamma_j^{-1}(i/n)\right\| \left | \bm{\nu}_j(i/n) \right |.
\end{align}
The rest of the proof follows from the exact reasoning as the arguments of the proof of Theorem 2.11 of \cite{DZar}. In particular, using (1)-(3) of Assumption \ref{assum_mainassumption}, by an argument similar to the equations between (D.31) and (D.32) of \cite{DZar}, we find that the first term of the right-hand side of (\ref{eq_normcontrol}) can be bounded by $\rO(j^{1.5}/n)$ and the second term can be bounded by $\rO(j^2/n).$ This completes our proof.

\end{Proof} 

\subsection{Proofs of Section \ref{sec_theorecticalanalysis}}

In this subsection, we prove the main results in Section \ref{sec_theorecticalanalysis}. 
\vspace*{4pt}

\begin{Proof}[\bf Proof of Theorem \ref{thm_consistency}] The proof is similar to that of Theorem 3.2 of \cite{ding2021simultaneous}. Due to similarity, we only sketch the key points. Similar to (\ref{eq_pacfestimate}), we set $\widehat{\bm{\phi}}_j(t)$ as the sieve estimator for $\bm{\phi}_j(t);$ that is
\begin{equation}\label{eq_olseestimate}
\widehat{\bm{\phi}}_j(t):=\widehat{\bm{\beta}}_j^* \mathbb{B}_j(t), \ \mathbb{B}_j(t)=\sum_{l=1}^j \mathbb{B}_{j,l}(t). 
\end{equation} 
Recall (\ref{eq_rhocxh}). For $1 \leq l \leq j,$ denote
\begin{equation*}
\bm{\phi}_{j,l}^{(c)}(t)=\sum_{k=1}^c a_{jk,l} \alpha_k(t), 
\end{equation*}
and $\bm{\phi}_j^{(c)}=\left(\bm{\phi}_{j,1}^{(c)},\cdots, \bm{\phi}_{j,j}^{(c)} \right)^* \in \mathbb{R}^j.$ Recall (\ref{eq_functionfinal}). The starting point is the following decomposition
\begin{equation*}
\| \widehat{\bm{\phi}}_j(t)- \bm{\phi}_{j,j}(t)\| \leq  \| \widehat{\bm{\phi}}_j(t)- \bm{\phi}^{(c)}_{j,j}(t)\|+|\bm{\phi}_{j,j}(t)-\bm{\phi}_{j,j}^{(c)}(t)|.  
\end{equation*}
Since the second term of the right-hand side of the above equation can be bounded by $\mathrm{O} \left(c^{-d} \right)$ using (\ref{eq_rhocxh}), it suffices to control the first term. According to (\ref{eq_olseestimate}), by Cauchy-Schwarz inequality, we see that 
\begin{equation}\label{eq_betadecomposition}
\| \widehat{\bm{\phi}}_{j,j}(t)- \bm{\phi}^{(c)}_{j,j}(t)\| \leq \sum_{l=1}^j \| \widehat{\bm{\phi}}_{j,l}(t)- \bm{\phi}^{(c)}_{j,l}(t)\| \leq \sqrt{j} \zeta_c  \| \bm{\beta}_j-\widehat{\bm{\beta}}_j \|. 
\end{equation} 
Moreover, for the OLS estimation, we have that
\begin{equation}\label{eq_eeeebbbbbbebbebebe}
 \bm{\beta}_j-\widehat{\bm{\beta}}_j=n(Y_j^* Y_j)^{-1} \frac{Y_j^* \bm{\epsilon}_j}{n}, 
\end{equation}
where $\bm{\epsilon}_j=(\epsilon_{j+1,j}+\mathfrak{r}, \cdots, \epsilon_{n,j}+\mathfrak{r})^* \in \mathbb{R}^{n-j}$ with $\mathfrak{r}$ representing the error term in (\ref{eq_representationxi}). The rest of the discussions follow lines of that of Theorem 3.2 of \cite{ding2021simultaneous}. 

To control the right-hand side of (\ref{eq_eeeebbbbbbebbebebe}), first, by an argument similar to (A.9) of \cite{ding2021simultaneous}, for $\Sigma^{(j)}$ defined in (\ref{eq_longruncovariancematrixdefinition}),  we find that
\begin{equation*}
\left\| \frac{1}{n} Y_j^* Y_j-\Sigma^{(j)} \right\|=\mathrm{O}_{\mathbb{P}} \left( jc \left( \frac{\xi_c^2}{\sqrt{n}}+\frac{\xi_c^2 n^{\frac{2}{\tau+1}}}{n} \right)  \right). 
\end{equation*} 
Combining the above discussion, together with Assumption \ref{assu_integratedeigenvalueassumption} and the assumption of (\ref{eq_rateassumption}), we see that 
\begin{equation}\label{eq_yjboundyjbound}
\left\| n(Y_j^* Y_j)^{-1} \right\|=\mathrm{O}_{\mathbb{P}}(1). 
\end{equation}
Second, according to a discussion similar to (A.12) of \cite{ding2021simultaneous} and the fact that $\mathfrak{r}=\mathrm{O}_{\mathbb{P}}(j^{2}/n+jc^{-d})$, we have that
\begin{equation}\label{eq_residualcontrol}
\left\| \frac{Y_j^* \bm{\epsilon}_j}{n} \right\|=\mathrm{O}_{\mathbb{P}} \left( \xi_c \left(1+ \frac{j^{2}}{n}+jc^{-d} \right)  \sqrt{\frac{jc}{n}} \right).
\end{equation}
We point out that the error rate $j^{2}/n$ is faster than the $j^{2.5}/n$ in (A.12) of \cite{ding2021simultaneous} or (2.23) and (2.24) of \cite{DZar} because of the assumption that $\{x_i\}$ is a mean-zero time series. Inserting the above two bounds into (\ref{eq_eeeebbbbbbebbebebe}), we conclude that 
\begin{equation*}
\| \bm{\beta}_j-\widehat{\bm{\beta}}_j \|=\mathrm{O}_{\mathbb{P}} \left( \xi_c \left(1+ \frac{j^{2}}{n}+jc^{-d} \right)  \sqrt{\frac{jc}{n}}\right).
\end{equation*}
Combining (\ref{eq_betadecomposition}), we can conclude our proof. 
\end{Proof}
\vspace*{4pt}

\begin{Proof}[\bf Proof of Theorem \ref{lem_controltbtrho0}]
The proof follows from lines of those of  Proposition 3.7 of \cite{DZar} verbatim. We omit the details. 
\end{Proof}
\vspace*{4pt}

\begin{Proof}[\bf Proof of Theorem \ref{lem_controltbtrho}]
The proof of part 1 follows from lines of those of  Proposition 3.7 of \cite{DZar} verbatim. We omit the details.

For part 2, we start with (\ref{eq_closenessstatistics}). Before proceeding to the actual proof, we first explore the relation between $\rho_k(t)$ and $\phi_{\mathsf{h},k}(t)$ for $1 \leq k \leq \mathsf{h}$ under the null hypothesis of (\ref{eq_testnotwo}). First of all, under (\ref{eq_testnotwo}), according to Theorem \ref{lem_pacf}, we can conclude that 
\begin{equation}\label{eq_mione}
\rho_{i,\mathsf{h}} \equiv \phi_{i\mathsf{h}, \mathsf{h}}=\rO\left( \min \left\{ \frac{\mathsf{h}^2}{n}, \left( \frac{\log \mathsf{h}+1}{\mathsf{h}} \right)^{\tau-1} \right\} \right). 
\end{equation} 
Combining with an argument similar to Theorem 2.11 of \cite{DZar} and the smoothness of $\phi_{\mathsf{h}, \mathsf{h}}(t)$, we can conclude that for all $0 \leq t \leq 1$
\begin{equation}\label{eq_mitwo}
\phi_{\mathsf{h}, \mathsf{h}}(t)=\rO \left( \frac{\mathsf{h}^2}{n} \right).
\end{equation}
Then with an argument similar to Theorem 2.11 of \cite{DZar}, in light of the representation (\ref{eq_jorderlocal}), we see that 
\begin{equation*}
x_i=\sum_{j=1}^{\mathsf{h}-1} \phi_{\mathsf{h},j}\left( \frac{i}{n} \right) x_{i-j}+\epsilon_{i,\mathsf{h}}+\rO_{\mathbb{P}} \left( \frac{\mathsf{h}^{2}}{n} \right).
\end{equation*}
Therefore, we can repeat (\ref{eq_mione}) and (\ref{eq_mitwo}) to get that 
\begin{equation}\label{eq_bbbbb}
\phi_{\mathsf{h},k}(t)=\rO\left( \frac{k^2}{n} \right), \ 1\leq k \leq \mathsf{h},
\end{equation}
and conclude that 
\begin{equation*}
x_i=\epsilon_{i,\mathsf{h}}+\rO_{\mathbb{P}} \left( \frac{\mathsf{h}^3}{n} \right).
\end{equation*}
Using our assumptions on $\mathsf{h}$, we find that the error term of the above approximation is always negligible with high probability. With the above preparation, we proceed to the proof. By definition, we have that 
\begin{equation*}
T_{\text{BP}}-T_2=\sum_{k=1}^{\mathsf{h}} \int_0^1  \left( \widehat{\rho}_k^2(t)-\widehat{\phi}_{\mathsf{h},k}^2(t) \right)  \mathrm{d}t.
\end{equation*}
By a discussion similar to Theorem  \ref{thm_consistency}, we have that (see Theorem 3.3 of \cite{DZ1}) for $1 \leq k \leq \mathsf{h}$
\begin{equation*}
\phi_{\mathsf{h},k}(t)=\widehat{\phi}_{\mathsf{h},k}(t)+\rO_{\mathbb{P}} \left( \Psi(\mathsf{h},c)  \right),
\end{equation*}
where we recall (\ref{eq_psiparameter}). Combining the above controls with Theorem \ref{thm_consistency}, under the null hypothesis (\ref{eq_testnotwo}), we see that
\begin{equation*}
\widehat{\rho}_k^2(t)=\rO_{\mathbb{P}} \left(\Psi^2(k,c) \right), \ \widehat{\phi}_{\mathsf{h},k}^2(t)=\rO_{\mathbb{P}} \left( \frac{k^4}{n^2}+\Psi^2(\mathsf{h},c) \right).
\end{equation*}
Consequently, we have that 
\begin{equation*}
T_{\text{BP}}-T_2=\mathrm{O}_{\mathbb{P}} \left( \mathsf{h} \Psi^2(\mathsf{h},c)+ \frac{\mathsf{h}^5}{n^2}+\sum_{k=1}^{\mathsf{h}} \Psi^2(k,c)\right). 
\end{equation*}
Then, for (\ref{eq_distribution}), using the definition of $g_k$ in (\ref{eq_defingk}) and the fact that $g_2 \asymp \sqrt{ \mathsf{h} c},$  the proof follows from  (\ref{eq_distribution11}),   (\ref{eq_closenessstatistics}) and the assumption of (\ref{eq_anotherassumption}).  

Finally, for part 3, by a discussion similar to (\ref{eq_mione}) and (\ref{eq_mitwo}), we find that (\ref{eq_alternativespecificform2ndapplication}) yields that
\begin{equation*}
\sum_{k=1}^{\mathsf{h}}\int_0^1 \phi_{\mathsf{h}, k}^2(t) \mathrm{d} t >C_{\alpha} \frac{\sqrt{\mathsf{h} c}}{n}.
\end{equation*}
Then the proof follows from lines of those of  Proposition 3.7 of \cite{DZar} verbatim.

\end{Proof}
\vspace*{4pt}

\begin{Proof}[\bf Proof of Corollary \ref{cor_boostrap}]
The proof follows directly from the proof of Theorem 3.10 of \cite{DZar} by setting $b_*=\mathsf{h}$ therein. 
\end{Proof}

\subsection{Additional assumptions}
In this section, we collect some more assumptions and provide more discussions on these technical assumptions.  The following assumption will be needed in our proof. It is a mild assumption and can be easily satisfied by many time series. We refer the readers to Section C.1 of \cite{DZar} for more details.
\vspace*{4pt}  

\begin{assu}\label{assu_ARbapproximation}
We assume the following assumptions hold true
\begin{enumerate}
\item Suppose  $\tau>4.5$ in Assumption \ref{assum_mainassumption}. Moreover, we assume that $c$ is of the form (\ref{eq_cform})  and satisfies that for large constant $C>0$
\begin{equation*}
\frac{C}{\tau}+a<1, \ \operatorname{and} \ da>2. 
\end{equation*}
\item  We assume that the derivatives of $\gamma(t,j)$ decay with $j$ as follows
\begin{equation*}
\sup_{t \in [0,1]} \sum_{j=0}^{\infty}|\gamma^{(d)}(t,j)|<\infty,
\end{equation*}
where $\gamma^{(d)}(t,j)$ is the d$th$ derivative of $\gamma(t,j)$ with respect to $t$.
\item   There exist constants $\omega_1, \omega_2 \geq 0,$ for some constant $C>0,$ we have 
\begin{equation*}
\sup_t | \nabla \mathbf{B}(t) | \leq C n^{\omega_1} c^{\omega_2}.
\end{equation*}
\end{enumerate}

\end{assu}

\section{Additional remarks and examples}\label{sec_additionalremark}

In this section, we provide several remarks and examples. The following remark provides more discussions on Theorem \ref{lem_pacf}.
\vspace*{4pt}
  
\begin{rem}
Theorem \ref{lem_pacf} is established for locally stationary time series as in Definition \ref{defn_locallystationary} where an exact cut-off is not available in general. An exception is the locally stationary AR(p) process as in \cite{zhou2013inference}, where
\begin{equation}\label{eq_timeseriesform}
x_i=\phi_{0}(i/n)+\sum_{j=1}^p \phi_j(i/n) x_{i-j}+\epsilon_i,
\end{equation} 
where $\{\epsilon_i\}$ is a time-varying white noise process. In this setting, it is clear that
\begin{equation*}
\rho_{i,j}=
\begin{cases}
\phi_j(i/n), & j=p ;\\
0, & j>p,
\end{cases}
\end{equation*}
so that $(\ref{eq_decayrate})$ and (\ref{eq_smoothapproximation}) holds trivially once $p$ is fixed or divergent slowly. In fact, as proved in Theorem 2.11 of \cite{DZar}, any locally stationary time series satisfying Definition \ref{defn_locallystationary} and Assumption \ref{assum_mainassumption} can be always well approximated by a time series in the form of (\ref{eq_timeseriesform}) with slowly diverging $p.$ In this regard, our results can be used to provide an order for the AR approximation; see Remark \ref{rem_app_determinep} below.     
\end{rem}

The following remark provides more explanations on the conditions in Assumption \ref{assum_mainassumption}. 
\vspace*{4pt}

\begin{rem}\label{rem_assumptionremark}
The conditions (1)--(3) in Assumption \ref{assum_mainassumption} are mild and commonly used in the literature. First, (1) is introduced to avoid the erratic behavior of the time series and frequently used in the statistics literature involving the covariance and precision matrix estimation \cite{cai2016,chen2013, DZ1, DZar,Yuan2010}. Moreover, as  proved in \cite[Proposition 2.9]{DZar},  it is equivalent to the uniform positiveness of the local spectral density function of $\{x_{i}\}.$ Second, (2) imposes the condition that the temporal structure of $\{x_{i}\}$ decays polynomially fast. This amounts to a short range dependent requirement for $\{x_{i}\}$  when $\tau>1.$ Analogous results can be easily obtained for the exponentially decay setting where  
\begin{equation*}
\max_{k,n}\left| \operatorname{Cov}(x_{k,n}, x_{k+r,n}) \right| \leq C a^{|r|}, \ \text{for some} \ 0<a<1. 
\end{equation*}  
Third, (3) requires that the  autocovariance functions of $\{x_{i}\}$ are smooth so that its PACFs can be estimated consistently. It is commonly used in the literature of locally stationary time series \cite{dahlhaus2012locally,dahlhaus2019towards,DZ1, DZar,WZ1}.   

\end{rem}

The remark below provides some insights on how to use our results to estimate the order of a locally stationary AR process. 
\vspace*{4pt}

\begin{rem}\label{rem_app_determinep}
 We discuss how to generalize the use of PACFs of stationary AR process to locally stationary AR process. For definiteness, we focus on the following time-varying AR($p$) process which has been used in \cite{zhou2013inference, dahlhaus2019towards}
\begin{equation}\label{eq_modellocalstationary}
x_i=\sum_{j=1}^p \phi_j(i/n) x_{i-j}+\epsilon_i, 
\end{equation} 
where $\{\epsilon_i\}$ is some locally-stationary white noise process and $\phi_j(\cdot), j=0,1,2,\cdots,p,$ are some smooth functions on $[0,1]$. Before estimating these time-varying coefficients, we first need to provide an estimator for $p.$ Here $p$ is allowed to diverge with $n.$ Based on the results of Theorem \ref{lem_controltbtrho}, inspired by the ideas in \cite{ding2022tracy,ding2023extreme},
we can propose a sequential test estimate based on the following hypothesis testing problem 
\begin{equation}\label{eq_testing4.1}
\mathbf{H}_0: p=p_0 \ \text{vs} \ \mathbf{H}_a: p_0<p \leq p_*, 
\end{equation} 
where $p_0$ is some pre-given integer representing our belief of the true value of $p$ and $p_*$ is some large integer that is interpreted as the maximum possible order the model can have. In light of (\ref{eq_firstnull}), we can use the following estimate 
\begin{equation}\label{eq_pestimateone}
\widetilde{p}=\max\{1 \leq j \leq p_*: \ \mathbf{H}_{a1} \ \text{in} \ (\ref{eq_firstnull}) \ \text{is accepted}\}. 
\end{equation}
\end{rem}

The following remark is related to the hypothesis testing (\ref{eq_testnotwo}).
\vspace*{4pt}

\begin{rem}
Two remarks are in order. First, in classic stationary time series analysis, in addition to Box-Pierce test, people also use the Ljung–Box (LB) test \cite{ljung1978measure}. Moreover, the BP and LB tests are asymptotic equivalent and follow the Chi-squared distribution with the same degree of freedom $\mathsf{h}$. In this regard, we can also modify  the LB test using  $$\mathcal{Q}_{\text{MLB}}:= n(n+2) \sum_{k=1}^{\mathsf h} \frac{\int_0^1  \widehat{\rho}_k(t)^2}{n-k}.$$ Moreover, we can study such a modified statistic as Theorem \ref{lem_controltbtrho}. Due to the asymptotic equivalence, we omit further details.  Second,  (\ref{eq_testnotwo}) is frequently used for model diagnostics. In this regard, it provides an alternative approach to choose the order of AR approximations by checking whether the residuals  follow white noise after fitting some AR models.
\end{rem}

Finally, we provide two frequently-used models of locally stationary time series in the literature and explain how Definition \ref{defn_locallystationary} and Assumption \ref{assum_mainassumption} can be easily satisfied.
\vspace*{4pt}

\begin{exam}\label{exam_timeseries}

We shall first consider the locally stationary time series model in \cite{WZ1,WZ2} using a physical representation so that 
\begin{equation}\label{eq_physcialrepresentation}
x_{i,n}=G_n(\frac{i}{n}, \mathcal{F}_i),\ i=1,2,\cdots, n, 
\end{equation}
where $\mathcal{F}_i=(\cdots, \eta_{i-1}, \eta_i)$ and $\eta_i, \ i  \in \mathbb{Z}$ are i.i.d centered random variables, and
$G_n:[0,1] \times \mathbb{R}^{\infty} \rightarrow \mathbb{R}$ is a measurable function such that $\xi_{i,n}(t):=G_n(t, \mathcal{F}_i)$ is a properly defined random variable for all $t \in [0,1].$  In (\ref{eq_physcialrepresentation}), by allowing the data generating mechanism $G_n$ depending on the time
index $t$ in such a way that $G_n(t,\mathcal{F}_i)$ changes smoothly with respect to $t$, one
has local stationarity in the sense that the subsequence $\{x_{i,n}
, . . . , x_{i+j-1,n} \}$ is
approximately stationary if its length $j$ is sufficiently small compared to $n$. Moreover, they quantify the temporal  decay using the physical dependence measure for (\ref{eq_physcialrepresentation}) as follows 
\begin{equation}\label{eq_dependencemeasure}
\delta(j,q):=\sup_{t \in [0,1]} || G_n(t,\mathcal{F}_0)-G_n(t, \mathcal{F}_{0,j}) ||_q.
\end{equation} 
Moreover, the following assumptions are needed to ensure local stationarity. 
\vspace*{4pt}

\begin{assu}\label{assum_local}
$G_n(\cdot, \cdot)$ defined in (\ref{eq_physcialrepresentation}) satisfies the property of stochastic Lipschitz continuity, i.e.,  for some $q>2$ and $C>0,$ 
\begin{equation}\label{assum_lip}
\left| \left| G_n(t_1, \mathcal{F}_{i})-G_n(t_2,\mathcal{F}_i) \right|\right|_q \leq C|t_1-t_2|, 
\end{equation} 
where $t_1, t_2 \in [0,1].$ Furthermore, 
\begin{equation}\label{assum_moment}
\sup_{t \in [0,1]} \max_{1 \leq i \leq n} ||G_n(t,\mathcal{F}_i) ||_q<\infty.
\end{equation} 
\end{assu}

It can be shown that time series $\{x_{i,n}\}$  with physical representation (\ref{eq_physcialrepresentation}) and Assumption \ref{assum_local} satisfies Definition \ref{defn_locallystationary}. In particular, for each fixed $t \in [0,1],$ $\gamma(t,j)$ in Definition \ref{defn_locallystationary} can be found easily using the following 
  \begin{equation}\label{eq_defncov}
\gamma(t,j)=\operatorname{Cov}(G_n(t, \mathcal{F}_0), G_n(t, \mathcal{F}_{j})).
\end{equation} 
Note that the assumptions
(\ref{assum_lip}) and (\ref{assum_moment}) ensure that   $\gamma(t,j)$ is Lipschiz continuous in $t$. Moreover, for each fixed $t,$ $\gamma(t,\cdot)$ is the autocovariance function of  $\{G_n(t,\cdot)\},$ which is a stationary process.

The physical representation form (\ref{eq_physcialrepresentation}) includes many commonly used locally stationary time series models. For example,  let $\{z_i\}$ be {zero-mean i.i.d. random variables (or a white noise) with variance $\sigma^2$.} We also assume $a_{j,n}(\cdot), j=0,1,\cdots$ be $C^d([0,1])$ functions such that
\begin{equation}\label{ex_linear}
G_n(t, \mathcal{F}_i)=\sum_{k=0}^{\infty} a_{k,n}(t) z_{i-k}. 
\end{equation} 
(\ref{ex_linear}) is a locally stationary linear process. It is easy to see that (2) and (3) of  Assumption \ref{assum_mainassumption}  will be satisfied if 
$ \sup_{t \in [0,1]} |a_{j,n}(t)| \leq C j^{-\tau}, \ j \geq 1; \ \sum_{j=0}^{\infty} \sup_{t \in [0,1]} |a_{j,n}'(t)|<\infty, $ and 
\begin{equation}\label{eq_coeffdecay}
 \sup_{t \in [0,1]} |a_{j,n}^{(d)}(t)|\leq  C j^{-\tau}, \ j \geq 1.
\end{equation} 

Furthermore, we note that the local spectral density function of \eqref{ex_linear} can be written as $f(t,w)= \sigma^2|\psi(t, e^{-\mathrm{i} j \omega})|^2,$ where $\psi(\cdot,\cdot)$ is defined such that $G_n(t, \mathcal{F}_i)=\psi(t, B) z_i$ with $B$ being the  backshift operator. As discussed in Remark \ref{rem_assumptionremark}, (1) of Assumption \ref{assum_moment}  will be satisfied if $|\psi(t, e^{-\mathrm{i} j \omega})|^2 \geq \kappa$ for all $t$ and $\omega,$ where $\kappa>0$ is some universal constant. For more examples of locally stationary time series in the form of (\ref{eq_physcialrepresentation}) especially nonlinear time series, we refer the readers to \cite{WW}, \cite[Section 2.1]{DZ1} , \cite[Example 2.2 and Proposition 4.4]{dahlhaus2019towards}, \cite[Proposition E.6]{karmakar2022simultaneous} and \cite{ding2021simultaneous, karmakar2022simultaneous,mayer2020functional}. Especially, the time-varying AR and ARCH models can be written into (\ref{ex_linear}) asymptotically \cite{DZar}, and  Assumptions \ref{assum_mainassumption} and \ref{assum_local} can be easily satisfied under mild assumptions. We refer the readers to the aforementioned references for more details.

For a second example,  note that in \cite{dahlhaus2006statistical,DPV,MR3097614}, the locally stationary time series is defined as follows (see Definition 2.1 of \cite{MR3097614}). $\{x_{i,n}\}$ is locally stationary time series if for each scaled time point $u \in [0,1],$ there exists a strictly stationary process $\{h_{i,n}(u)\}$ such that 
\begin{equation}\label{eq_definition11111}
|x_{i,n}-h_{i,n}(u)| \leq \left(|t_i-u|+\frac{1}{n} \right)U_{i,n}(u),  \ \text{a.s},
\end{equation} 
where $U_{i,n}(u) \in L^q([0,1])$ for some $q>0.$ By similar arguments as those of model \eqref{eq_physcialrepresentation} \cite{DZar},  Definition \ref{defn_locallystationary} as well as assumptions of this subsection can be verified for \eqref{eq_definition11111}, especially (\ref{eq_definition11111}) implies (\ref{eq_covdefn}).
\end{exam}

%
%
%
%
%

\section{Tuning parameters selection} \label{sec:choiceparameter}
In this section, we explain how to choose the tuning parameters associated with our proposed methodology.

First, we discuss how to choose the tuning parameters $c$ and $m$ used in Algorithm \ref{alg:boostrapping}.   We use a data-driven procedure proposed in \cite{bishop2013pattern} to choose $c.$  For a given integer $l,$ say $l=\lfloor 3 \log_2 n \rfloor,$ we divide the time series into two parts: the training part $\{x_i\}_{i=1}^{n-l}$ and the validation part $\{x_i\}_{i=n-l+1}^n.$  With some preliminary initial value $c$, we propose a sequence of candidates  $ c_j, \ j=1,2,\cdots, v,$ in an appropriate neighborhood of $c$ where $v$ is some given integer. For each of the choices $c_j,$  we fit a time-varying AR($h_2$) model as in (\ref{eq_residualestimate}) with $c_j$ sieve basis expansion using the training data set.  Then using the fitted model, we forecast the time series in the validation part of the time series.  Let $\widehat x_{n-l+1,j}, \cdots, \widehat x_{n,j}$ be the forecast of $x_{n-l+1},..., x_n,$ respectively using the parameter $c_j$. Then we choose the parameter $c_{j_0}$ with the minimum sample MSE of forecast, i.e.,
 \begin{equation*} 
{j_0}:= \argmin_{(j: 1 \leq j \leq v)} \frac{1}{l}\sum_{k=n-l+1}^n (x_k-\widehat x_{k,j})^2.
 \end{equation*}

{
To choose an $m$ for practical implementation, in \cite{ZZ1}, the author used the minimum volatility (MV) method to choose the window size $m$ for the scalar covariance function. The MV method does not depend on the specific form of the underlying time series dependence structure and hence is robust to misspecification of
the latter structure \cite{politis1999subsampling}. The MV method utilizes the fact that the covariance structure of $\widehat{\Pi}$ becomes stable when the
block size $m$ is in an appropriate range, where $\widehat{\Pi}=E[\widehat{\Phi} \widehat{\Phi}^*|(x_1,\cdots,x_n)]$ is defined as 
{ 
\begin{equation}\label{eq_widehatomega}
\widehat{\Pi}:=\frac{1}{(n-m-h_2+1)m} \sum_{i=h_2+1}^{n-m} \Big[ \Big(\sum_{j=i}^{i+m} \widehat{\bm{w}}_{h_2,j} \Big) \otimes \Big( \mathbf{B}(\frac{i}{n}) \Big) \Big] \times \Big[ \Big(\sum_{j=i}^{i+m} \widehat{\bm{w}}_{h_2,j} \Big) \otimes \Big( \mathbf{B}(\frac{i}{n}) \Big) \Big]^*.
\end{equation}
}    Therefore, it desires to minimize the standard errors of the latter covariance structure in a suitable range of candidate $m$'s. In detail, for a give large value $m_{n_0}$ and a neighborhood control parameter $h_0>0,$  we can choose a sequence of window sizes $m_{-h_0+1}<\cdots<m_1< m_2<\cdots<m_{n_0}<\cdots<m_{n_0+h_0}$  and obtain $\widehat{\Pi}_{m_j}$ by replacing $m$ with $m_j$ in (\ref{eq_widehatomega}), $j=-h_0+1,2, \cdots, n_0+h_0.$ For each $m_j, j=1,2,\cdots, m_{n_0},$ we calculate the matrix norm error of $\widehat{\Omega}_{m_j}$ in the $h_0$-neighborhood, i.e., 
\begin{equation*}
\mathsf{se}(m_j):=\mathsf{se}(\{ \widehat{\Pi}_{m_{j+k}}\}_{k=-h_0}^{h_0})=\left[\frac{1}{2h_0} \sum_{k=-h_0}^{h_0} \| \overline{\widehat{\Pi}}_{m_j}-\widehat{\Pi}_{m_j+k} \|^2 \right]^{1/2},
\end{equation*}
where $\overline{\widehat{\Pi}}_{m_j}=\sum_{k=-h_0}^{h_0} \widehat{\Pi}_{m_j+k} /(2h_0+1).$
Therefore, we choose the estimate of $m$ using 
\begin{equation*}
\widehat{m}:=\argmin_{m_1 \leq m \leq m_{n_0}} \mathsf{se}(m).
\end{equation*}
Note that in \cite{ZZ1} the author used $h_0=3$ and we also adopt this choice in the current paper. 
}

{Second, we discuss how to choose  a large value of $\mathsf{h}$ to construct the statistic $T_2$ in (\ref{eq_defnt2}). Theoretically, the lower bound for the order of $\mathsf{h}$ is given by $j^*$ as in (\ref{j_assumption}), while the upper bound is provided in the assumption of (\ref{eq_anotherbound}).  The lower bound yields that for all $j \geq \mathsf{h},$ $\sup_t \rho_j(t)=\mathrm{o}(n^{-1/2}),$ while the upper bound guarantees that the error term is negligible. To balance these two conditions, for practical implementation, we use the following value
\begin{equation*}
\mathsf{h}:=\min\left\{ 1 \leq j \leq \mathsf{h}^*: \ \mathbf{H}_{01} \ \text{in} \ (\ref{eq_firstnull}) \ \text{is accepted}  \right\},
\end{equation*}
where $\mathsf{h}^*>0$ is some pre-given large value (say, $\mathsf{h}^*=50$). We emphasize that this can be easily done using our $\mathtt{R}$ package $\mathtt{Sie2nts}$ by generating a plot as in Figure \ref{fig_pvalueplot}. 
}

\section{Additional simulation results}\label{sec_additionalsimulationresults} 
In this section, we provide additional numerical simulation results.  
\subsection{More results on other types of models}\label{sec_simulationmore}
In this section, we conduct more numerical simulations  using both stationary and non-stationary MA(1) models.  For some constant $\delta \in [0,0.5],$ we consider the stationary MA(1) process 
\begin{equation}\label{eq_mastationary}
x_i=\epsilon_i+\delta \epsilon_{i-1},
\end{equation}
and the locally stationary MA(1) process 
\begin{equation}\label{eq_manonstationary}
x_i=\epsilon_i+\delta \sin(2 \pi i/n) \epsilon_{i-1},
\end{equation}
where $\epsilon_i, 1 \leq i \leq n,$ are i.i.d. standard Gaussian random variables. Note that when $\delta=0$ in (\ref{eq_mastationary}) and (\ref{eq_manonstationary}), they both reduce to the standard white noise. Since (\ref{eq_mastationary}) and (\ref{eq_manonstationary}) are essentially AR($\infty$) models, one can follow the discussions of Section 3.3 of \cite{shumway2017time} to calculate the true PACFs.

\begin{figure}[!ht]
\hspace*{0cm}
\begin{subfigure}{0.55\textwidth}
\includegraphics[width=6.5cm,height=6cm]{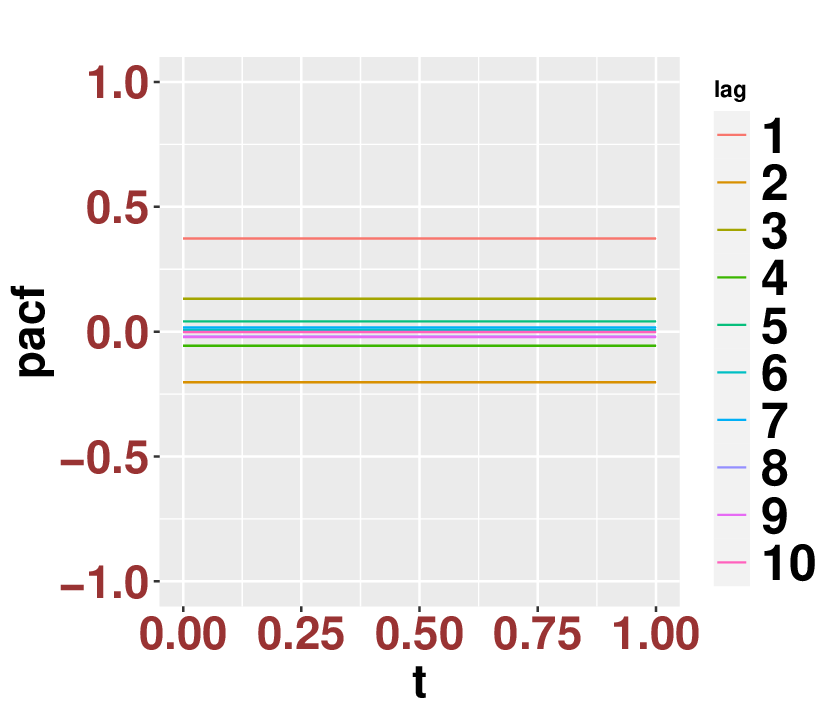}
\caption{PACFs for model (\ref{eq_mastationary}).}\label{subfig_nullorthogonaltypei200}
\end{subfigure}
\hspace{0cm}
\begin{subfigure}{0.55\textwidth}
\includegraphics[width=6.5cm,height=6cm]{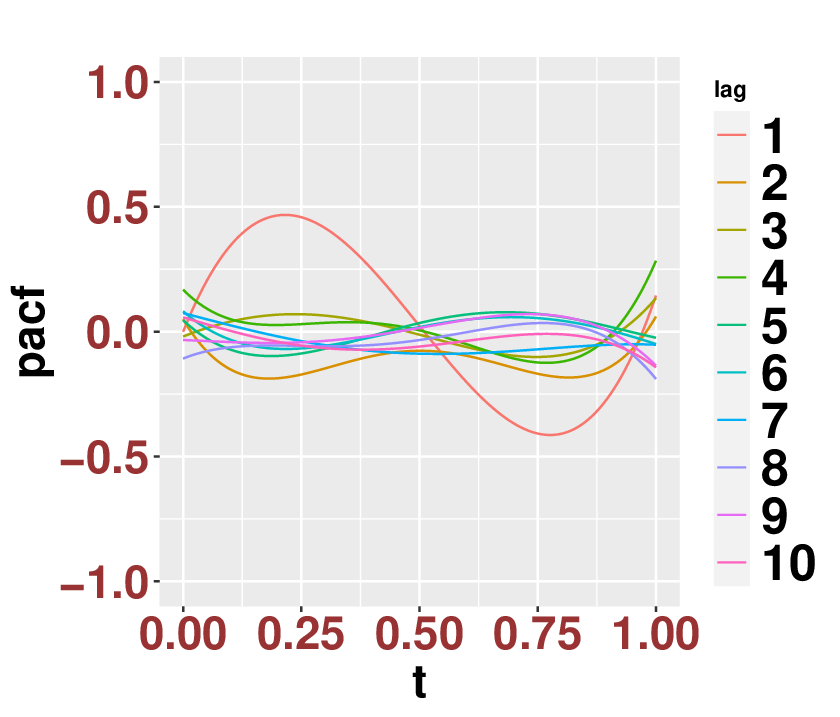}
\caption{PACFs for model (\ref{eq_manonstationary}).}\label{subfig_nullorthogonaltypei500}
\end{subfigure}
\vspace*{-0.8cm}
\caption{{ PACF plots for models (\ref{eq_mastationary}) and (\ref{eq_manonstationary}). Here $n=600.$  }}
\label{fig_pacfplotma}
\end{figure}

\begin{figure}[!ht]
\hspace*{0cm}
\begin{subfigure}{0.55\textwidth}
\includegraphics[width=6.5cm,height=6cm]{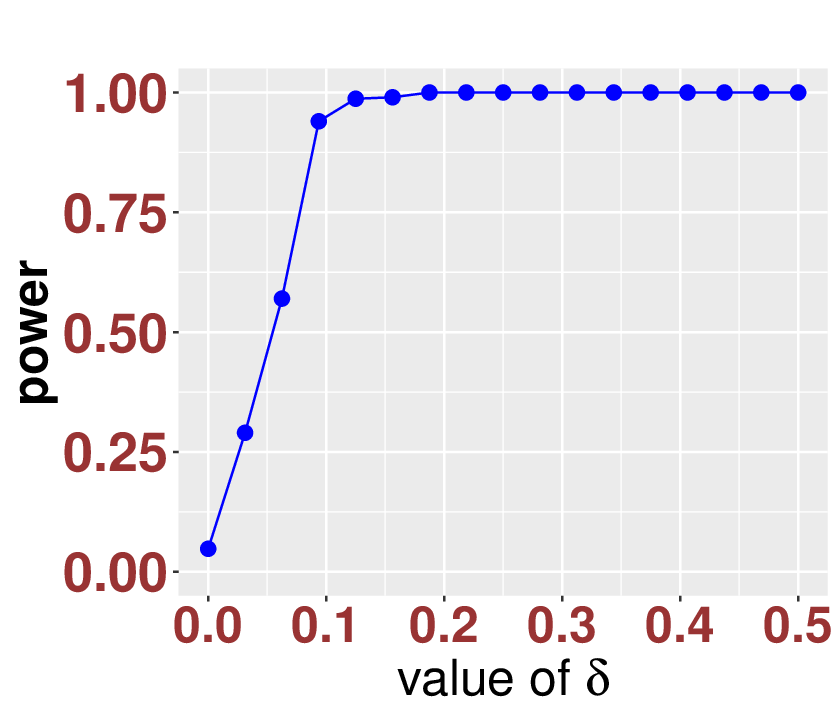}
\caption{Power for model (\ref{eq_mastationary}).}\label{subfig_nullorthogonaltypei200}
\end{subfigure}
\hspace{0cm}
\begin{subfigure}{0.55\textwidth}
\includegraphics[width=6.5cm,height=6cm]{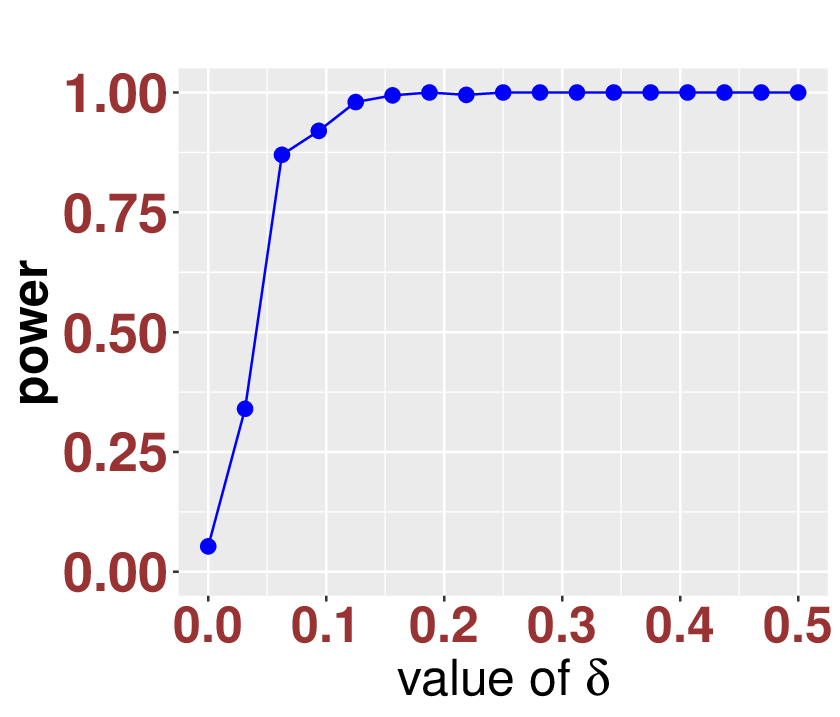}
\caption{Power for model (\ref{eq_manonstationary}).}\label{subfig_nullorthogonaltypei500}
\end{subfigure}
\vspace*{-0.8cm}
\caption{{ Power for models (\ref{eq_mastationary}) and (\ref{eq_manonstationary}) under the alternative of (\ref{eq_nullandalternativerespectively}). Here the type I error rate $\alpha=0.05$ and $n=600$. We use the Legendre polynomials as the basis functions. The results are reported based 1,000 repetitions.}}
\label{fig_powerplotma}
\end{figure}

For the estimation of the PACFs, in Figure \ref{fig_pacfplotma}, we provide the plots of the PACFs of the first 10 lags for 
both (\ref{eq_mastationary}) and (\ref{eq_manonstationary}). We can see that our estimates are reasonably accurate. Regarding the inference of the PACFs, for the purpose of definiteness, we focus on the white noise test (\ref{eq_nullandalternativerespectively}) where $\mathbf{H}_0$ corresponds to $\delta=0$ in (\ref{eq_mastationary}) and (\ref{eq_manonstationary}) and $\mathbf{H}_a$ corresponds to an MA(1) alternative that $\delta>0.$ Under $\mathbf{H}_0,$ (\ref{eq_mastationary}) and (\ref{eq_manonstationary}) are essentially the same model. For $n=600,$ under the type I error rate $\alpha=0.05,$ the simulated type I error rates are $0.047, 0.051, 0.056$ for the Fourier, Legendre and Daubechies-9 basis functions, respectively based on 1,000 repetitions. This shows the accuracy of our test. To examine the power, in Figure \ref{fig_powerplotma}, we report how the simulated power changes when $\delta$ deviates away from zero. We can conclude that our proposed test is reasonably powerful once the alternative deviates from the null.

\subsection{Comparison with \cite{10.1214/20-EJS1748} on estimating the PACFs}\label{sec_simulationcomparison}

In this section, we compare our method with the ones proposed in \cite{10.1214/20-EJS1748} in terms of the estimation of the PACFs using the mean integrated squared error (MISE). To implement \cite{10.1214/20-EJS1748}, we use the $\mathtt{R}$ package $\mathtt{lpacf}$ which is developed by the authors of  \cite{10.1214/20-EJS1748}.

For definiteness, we consider the AR type models (\ref{eq_stationarymodel}) and (\ref{eq_nonstationarymodel}) with $\delta_1=0.5$ and $\delta_2=0.$ Consequently, for model (\ref{eq_stationarymodel}), the true PACFs are 
\begin{equation*}
\rho_j(t) \equiv 
\begin{cases}
0.5 & j=1 \\
0 & j \geq 2 
\end{cases},
\end{equation*}
and for model (\ref{eq_nonstationarymodel}), the true PACFs are 
\begin{equation*}
\rho_j(t)=
\begin{cases}
0.5 \sin (2 \pi t) & j=1 \\
0 & j \geq 2
\end{cases}.
\end{equation*}

Let $\widehat{\rho}_j(t)$ be some estimator for $\rho_j(t),$ the MISE is defined as 
\begin{equation*}
\mathsf{MISE}(j)=\int_0^1 (\widehat{\rho}_j(t)-\rho_j(t))^2 \mathrm{d} t.
\end{equation*}
In the actual calculations, MISE can be well approximated by an Riemann summation.
In Table \ref{table_comparision}, we record the $\mathsf{MISE}(j), j=1,2,3,4,$ for our proposed method (denoted as Proposed) and the two methods in \cite{10.1214/20-EJS1748} (the wavelet-based method is denoted as Lpacf-I and  the Epanechnikov windowed method is denoted as Lpacf-II).  We can conclude that our proposed method has better finite sample performance than \cite{10.1214/20-EJS1748} which are known to have worse performance near the boundaries.

\begin{table}[!ht]
\begin{center}
\setlength\arrayrulewidth{1pt}
\renewcommand{\arraystretch}{1.5}
{\fontsize{10}{11}\selectfont 
\begin{tabular}{|c|cccc|}
\hline
Methods/Lags & \multicolumn{1}{c|}{$j=1$} & \multicolumn{1}{c|}{$j=2$} & \multicolumn{1}{c|}{$j=3$} & \multicolumn{1}{c|}{$j=4$}  \\ 
\hline
     & \multicolumn{4}{c|}{Model (\ref{eq_stationarymodel})}                                                                                                                                                                                                                                                                                          \\
   \hline
Proposed     &          $8 \times 10^{-4}$  &  $7 \times 10^{-4}$ &    $3 \times 10^{-3}$                   &         $1 \times 10^{-3}$                 \\
Lpacf-I    &     $9 \times 10^{-3}$       & $7 \times 10^{-3}$                          &                          $9 \times 10^{-3}$  &   $7 \times 10^{-3}$                                                \\
Lpacf-II    &  0.042 & 0.054   & 0.043                        &      0.051              \\
\hline
      & \multicolumn{4}{c|}{Model (\ref{eq_nonstationarymodel})}                                                                                                                                                                                                                                                                                         \\
       \hline
Proposed    &         $9.8 \times 10^{-3}$    & $1.9 \times 10^{-3}$ &           $2.5 \times 10^{-3}$ &      $2.7 \times 10^{-3}$                        \\
Lpacf-I    &    0.016    & $7.5  \times 10^{-3}$    &                  $3.5  \times 10^{-3}$                             &          $4  \times 10^{-3}$                 \\
Lpacf-II    &  0.049  & 0.051    &      0.042                   &                         0.055   \\
 \hline
\end{tabular}
}
\end{center}
\vspace*{-0.4cm}
\caption{Comparison of the accuracy of the estimation of  PACFs using MISE. Here Proposed stands for our estimator (\ref{eq_pacfestimate}) and Lpacf-I and Lpacf-II are the two methods introduced in \cite{10.1214/20-EJS1748}. Our proposed method can be implemented using our package $\mathtt{Sie2nts}$ and Lpacf-I/II can be implemented using the package $\mathtt{lpacf}.$ Here $n=1024.$    
}
\label{table_comparision}
\end{table}


 \vspace*{-0.65cm}

\bibliographystyle{biometrika}
\bibliography{PACF,corrtest}

\end{document}